\documentclass[11pt,reqno]{amsart}

\usepackage{amsthm}
\usepackage{amssymb}
\usepackage{latexsym}
\usepackage{multicol}
\usepackage{verbatim,enumerate}
\usepackage{hyperref}
\usepackage{nicematrix}
\usepackage{amsmath, amscd}
\usepackage{mathrsfs}
\usepackage[all,cmtip]{xy}
\usepackage{soul}

\advance\textwidth by 1.2in \advance\oddsidemargin by -.6in \advance\evensidemargin by -.6in
\usepackage{tikz}
\usetikzlibrary{decorations.markings}
\tikzstyle{vertex}=[ circle, fill, draw, inner sep=0pt, minimum size=4pt,]
\tikzstyle{edge}= [thick]

\newtheorem*{cor}{Corollary}
\newtheorem*{lem}{Lemma}
\newtheorem*{prop}{Proposition}

\theoremstyle{definition} 
\theoremstyle{definition}
\newtheorem{thm}{Theorem}
\newtheorem*{thm*}{Theorem}

\newtheorem*{rem}{Remark}

\newenvironment{pf}{\proof}{\endproof}
\newcounter{cnt}
\newenvironment{enumerit}{\begin{list}{{\hfill\rm(\roman{cnt})\hfill}}{%
\settowidth{\labelwidth}{{\rm(iv)}}\leftmargin=\labelwidth%
\advance\leftmargin by \labelsep\rightmargin=0pt\usecounter{cnt}}}{\end{list}} \makeatletter
\def\mydggeometry{\makeatletter\dg@YGRID=1\dg@XGRID=20\unitlength=0.003pt\makeatother}
\makeatother \theoremstyle{remark}

\numberwithin{equation}{section}

\newcommand{\wt}{\operatorname{wt}}

\newcommand{\nc}{\newcommand}
\newcommand{\rnc}{\renewcommand}

\nc{\cal}{\mathcal} \nc{\goth}{\mathfrak} \rnc{\bold}{\mathbf}

\renewcommand{\Bbb}{\mathbb}
\nc\bomega{{\mbox{\boldmath $\omega$}}} \nc\bpsi{{\mbox{\boldmath $\Psi$}}}
 \nc\sing{{\rm sing}}
 \nc\balpha{{\mbox{\boldmath $\alpha$}}}
 \nc\bbeta{{\mbox{\boldmath $\beta$}}}
 \nc\bpi{{\mbox{\boldmath $\pi$}}}
  \nc\bpis{{\mbox{\boldmath \scriptsize$\pi$}}}
 \nc\bullets{{\mbox{\scriptsize $\bullet$}}}
 \nc\bvarpis{{\mbox{\boldmath \scriptsize$\varpi$}}}
  \nc\bvarpi{{\mbox{\boldmath $\varpi$}}}

\nc\bepsilon{{\mbox{\boldmath $\epsilon$}}}

  \nc\bomegas{{\mbox{\boldmath\scriptsize $\omega$}}}
  \nc\bepsilons{{\mbox{\boldmath \scriptsize$\epsilon$}}}
\nc{\spi}{{\rm sp}}
\nc\hlien{\hat{\lie n}^+}
  \nc\btaus{{\mbox{\boldmath \scriptsize$\tau$}}}\nc\bxi{{\mbox{\boldmath $\xi$}}}
\nc\bmu{{\mbox{\boldmath $\mu$}}} \nc\bcN{{\mbox{\boldmath $\cal{N}$}}} \nc\bcm{{\mbox{\boldmath $\cal{M}$}}} \nc\blambda{{\mbox{\boldmath
$\lambda$}}}
\nc\btau{{\mbox{\boldmath
$\tau$}}}

\newcommand{\lie}[1]{\mathfrak{#1}}

\makeatletter
\def\section{\def\@secnumfont{\mdseries}\@startsection{section}{1}%
  \z@{.7\linespacing\@plus\linespacing}{.5\linespacing}%
  {\normalfont\scshape\centering}}
\def\subsection{\def\@secnumfont{\bfseries}\@startsection{subsection}{2}%
  {\parindent}{.5\linespacing\@plus.7\linespacing}{-.5em}%
  {\normalfont\bfseries}}
\makeatother

 \nc{\Hom}{\operatorname{Hom}}
  \nc{\mode}{\operatorname{mod}}
\nc{\End}{\operatorname{End}} \nc{\wh}[1]{\widehat{#1}} \nc{\Ext}{\operatorname{Ext}}
 \nc{\ch}{\operatorname{ch}} \nc{\ev}{\operatorname{ev}}
\nc{\Ob}{\operatorname{Ob}} \nc{\soc}{\operatorname{soc}} \nc{\rad}{\operatorname{rad}} \nc{\head}{\operatorname{head}}

 \nc{\Cal}{\cal} \nc{\Xp}[1]{X^+(#1)} \nc{\Xm}[1]{X^-(#1)}
\nc{\on}{\operatorname} \nc{\Z}{{\bold Z}} \nc{\J}{{\cal J}} \nc{\C}{{\bold C}} \nc{\Q}{{\bold Q}}

\nc{\N}{{\Bbb N}} \nc\boa{\bold a} \nc\bob{\bold b} \nc\boc{\bold c} \nc\bod{\bold d} \nc\boe{\bold e} \nc\bof{\bold f} \nc\bog{\bold g}
\nc\boh{\bold h} \nc\boi{\bold i} \nc\boj{\bold j} \nc\bok{\bold k} \nc\bol{\bold l} \nc\bom{\bold m} \nc\bon{\bold n} \nc\boo{\bold o}
\nc\bop{\bold p} \nc\boq{\bold q} \nc\bor{\bold r} \nc\bos{\bold s} \nc\boT{\bold t} \nc\boF{\bold F} \nc\bou{\bold u} \nc\bov{\bold v}
\nc\bow{\bold w} \nc\boz{\bold z} \nc\boy{\bold y} \nc\ba{\bold A} \nc\bb{\bold B} \nc\bc{\mathbb C} \nc\bd{\bold D} \nc\be{\bold E} \nc\bg{\bold
G} \nc\bh{\bold H} \nc\bi{\bold I} \nc\bj{\bold J} \nc\bk{\bold K} \nc\bl{\bold L} \nc\bm{\bold M}  \nc\bo{\bold O} \nc\bp{\bold
P} \nc\bq{\bold Q} \nc\br{\bold R} \nc\bs{\bold S} \nc\bt{\bold T} \nc\bu{\bold U} \nc\bv{\bold V} \nc\bw{\bold W} \nc\bx{\bold
x} \nc\KR{\bold{KR}} \nc\rk{\bold{rk}} \nc\het{\text{ht }}
\nc\bz{\mathbb Z}
\nc\bn{\mathbb N}
\nc\us{\underline \bos}
\nc\uS{ \bs_{\rm alt}}
\nc\pr{\rm pr}

\nc\toa{\tilde a} \nc\tob{\tilde b} \nc\toc{\tilde c} \nc\tod{\tilde d} \nc\toe{\tilde e} \nc\tof{\tilde f} \nc\tog{\tilde g} \nc\toh{\tilde h}
\nc\toi{\tilde i} \nc\toj{\tilde j} \nc\tok{\tilde k} \nc\tol{\tilde l} \nc\tom{\tilde m} \nc\ton{\tilde n} \nc\too{\tilde o} \nc\toq{\tilde q}
\nc\tor{\tilde r} \nc\tos{\tilde s} \nc\toT{\tilde t} \nc\tou{\tilde u} \nc\tov{\tilde v} \nc\tow{\tilde w} \nc\toz{\tilde z} \nc\woi{w_{\omega_i}}
\nc\chara{\operatorname{Char}}

\begin{document}

\title[Dominant $\ell$--weights and  maps between Weyl modules for $\widehat\bu_q(\lie{sl}_{n+1})$]{On dominant $\ell$--weights and  maps between Weyl modules for quantum affine $A_n$}
\author{Matheus Brito}
\address{Departamento de Matematica, UFPR, Curitiba - PR - Brazil, 81530-015}
\email{mbrito@ufpr.br}
\thanks{M.B. was partially supported CNPq grant 405793/2023-5}
\author{Vyjayanthi Chari}
\address{Department of Mathematics, University of California, Riverside, 900 University Ave., Riverside, CA 92521, USA}
\email{chari@math.ucr.edu}
\thanks{V.C. was partially supported by a travel grant from the Simons Foundation.}

\begin{abstract} We  determine the set of dominant $\ell$--weights in the Weyl (or standard) modules for quantum affine $A_n$. We then prove that the space of homomorphisms between standard modules is at most one-dimensional and give a necessary and sufficient condition for equality to hold. We also describe the socle of the standard module and prove that the socle is simple for large $n$. Finally, we give applications of our results to mixed Weyl modules, calculating extensions in the category and identify new families of tensor subcategories of finite dimensional representations.
    
\end{abstract}

\maketitle

\section*{Introduction}

The category of finite-dimensional representations of a quantum affine algebra $\widehat\bu_n$ (associated to a simple Lie algebra $\lie g$) has been studied intensively for a long time. There are still many interesting open questions about the structure of the representations
  and the homological properties of the category are largely unexplored. The goal of this paper is to address some of the basic questions in the case when  $\lie g$ is of type $A_n$. \\\\
  It was shown in \cite{CP01} that there exists a universal family of finite-dimensional modules, called Weyl  modules which are  indexed by the same set as the irreducible objects  of the category.  In many ways, the Weyl modules, although finite-dimensional, play a role analogous to the Verma modules in the Bernstein--Gelfand--Gelfand category $\cal O(\lie g)$. For instance, the results of \cite{C01}, \cite{VV02} show that the Weyl  modules admit a nice tensor product decomposition. This allows one to compute their $q$--characters in an explicit way.
  The irreducible modules on the other hand have complicated structure in general and their characters are alternating triangular sums of the character of the Weyl modules. \\\\
  Proving even these results, in the quantum affine case, is more complicated than in the case of $\cal O$ where the character of the Verma module is easily computed. This is because the Verma modules are free over a nice subalgebra; the Weyl modules by contrast have no such property. This difference also causes difficulties in other natural questions that one could ask about Weyl modules. For instance, one of the starting points in the study of the category $\cal O$ is the understanding of space of maps between two Verma modules and also proving the existence of a simple socle.
  In this paper we study the analogs of these results for Weyl modules in type $A_n$. We restrict our attention to  the Hernandez--Leclerc subcategory $\mathscr F_n$ of finite-dimensional modules. This subcategory arises naturally in connections to monoidal categorification of cluster algebras and quiver Hecke algebras (\cite{HL10,HL13a,KKOP,KKOP22}) and exhibits a rich combinatorial structure. They are also connected through a Schur--Weyl type duality with  complex smooth representations of $GL_N(F)$, where $F$ is a non-Archimedean field (\cite{Gur21,LM14,LM18}). \\\\
  Our main results are the following. We prove that $\dim\Hom_{\widehat\bu_n}(W_1,W_2)\le 1$ where $W_1$ and $W_2$ are Weyl modules and we give a necessary and sufficient condition for equality to hold. We show  that any non-zero between Weyl modules is injective. We  give an algorithm to determine the socle and  prove that the irreducible modules which occur in the socle are all Weyl modules and occur multiplicity one. Under suitable conditions we prove that the socle is simple and describe it explicitly. We emphasize here that the socle is not always a fundamental Weyl module; this is in sharp contrast to what happens in the classical limit, see for instance \cite{CL06}, \cite{FoL07}.\\\\
  In the final section of the paper we give some applications of our results. We show that one can define canonically, families of tensor subcategories of $\mathscr F_n$, generalizing the category $\mathscr C_\ell$ of \cite{HL10}. We also  study  mixed Weyl modules; these are modules which have the same character as the Weyl modules but are not isomorphic to a Weyl module. We identify certain canonical modules which occur in the head and the socle of a Weyl module.  The work of \cite{GM25} shows that  even in the case of $A_1$ the socle of a mixed Weyl module will have multiplicity. 
Finally, we  give a sufficient condition for $\Ext_{\widehat\bu_n}^1(V,W)$ to be trivial  where $V,W$ are  modules in $\mathscr F_n$.\\\\
The extension of our results to the  quantum affine algebra associated to other types is  complicated; in part this is because the structure and character of the fundamental Weyl modules is more complicated and the theory of  $q$--characters for the irreducible modules is not as well--developed. We hope to return to these questions in  subsequent papers.

\section{Main Results}
We  introduce the definitions of the main objects of study and then state the main results of the paper.
\subsection{}  As usual $\mathbb C$ (resp. $\mathbb C^\times$, $\mathbb Z$, $\mathbb Z_+$, $\mathbb N$) will denote the set of complex numbers (resp. non-zero complex numbers, integers, non-negative integers, positive integers).  Given any $r\in\mathbb N$ we let $\Sigma_r$ be the group of permutations on $r$ letters and $\sigma_p$, $1\le p\le r-1$ denote the simple transpositions. \\\\
Assume throughout that $q$ is a non--zero complex number and not a root of unity.
\subsection{The algebra $\widehat\bu_n$ and the category $\mathscr F_n$}\label{basicdef}  For $n\in\mathbb N$, let  $\widehat\bu_n$ be  the quantum loop algebra associated to $\widehat{\lie{sl}}_{n+1}(\mathbb C)$;  we refer the reader to \cite{CP91,CP95} for precise definitions. For our purposes, it is enough to recall that $\widehat\bu_n$ is a Hopf algebra with 
an infinite set of generators: $x_{i,s}^\pm$, $k_i^{\pm 1}$, $\phi^\pm_{i,s}$,  $1\le i\le n$ and $s\in\mathbb Z$. We let $\widehat\bu_n^0$ be the subalgebra of $\widehat\bu_n$ generated by $\phi_{i,s}^\pm$, $1\le i\leq n$ and $s\in\mathbb Z$. 
\\\\
It is well known (see \cite{CP91,CP95}) that the irreducible finite-dimensional modules of $\widehat\bu_n$  are parameterized  by elements of a free abelian monoid with identity $\bold 1$   and generators  $\bvarpi_{m,a}$, $1\le m\le n$ and $a\in\mathbb C^{\times}$. The trivial representation of $\widehat\bu_n$ corresponds to the identity element of the monoid. Moreover, it was shown in \cite{CP01} that corresponding to an element of this monoid there also exists a finite-dimensional indecomposable module called a Weyl module which has the irreducible module as its unique irreducible quotient. 
\\\\
Let $\mathscr F_{n}$ be the full subcategory of the category of finite-dimensional representations of $\widehat\bu_n$ consisting of objects whose  Jordan--Holder components are indexed by elements of the form $\bvarpi_{m_1,q^{a_1}}\cdots\bvarpi_{m_r,q^{a_r}}$ with $a_s-m_s\in2\mathbb Z$ for $1\le s\le r$. It was proved in \cite{HL10} that $\mathscr F_{n}$ is a rigid tensor category and  we  let $\cal K_0(\mathscr F_{n})$ be the corresponding  Grothendieck ring and for  $V\in \mathscr F_n$ we denote by $[V]\in \cal K_0(\mathscr F_n)$. The results of \cite{FR99} show that this     ring is commutative with basis given by the classes of the simple objects.
\subsection{The group $\cal I_n$} It will be convenient for us to use a different index set for the simple objects of $\mathscr F_n$. In fact, here we borrow the language of representations of $GL_n(F)$ where $F$ is a non--Archimedean field (see \cite{BC25} and \cite{Gur21} for further details).

Let  $\mathbb I_n$ be the set of intervals $[i,j]$ with $i,j\in\mathbb Z$ and $0\le j-i\le n+1$ and for $r\ge 1$ let $\mathbb I_n^r$ be the set of ordered $r$--tuples of elements of $\mathbb I_n$. Define $\cal I_n^+$ (resp. $\cal I_n$) to  be the  free abelian monoid (group) with identity $\bold 1$ and generators
 $\bomega_{i,j}$ with $[i,j]\in\mathbb I_n$, $1\leq j-i\leq n$. We understand that $\bomega_{i,i}=\bomega_{i,i+n+1}=\bold 1$ for all $i\in \mathbb Z$.
 We have a map $\mathbb I_n^r\to\cal I_n^+$ given by  $$ \bos=([i_1,j_1],\cdots, [i_r,j_r])\mapsto \bomega_{\bos}=\bomega_{i_1,j_1}\cdots\bomega_{i_r, j_r}.$$
Identifying a  pair $(m,q^a)$ satisfying  $a-m\in2\mathbb Z$ with the interval $[\frac12(a-m), \frac12(a+m)]$ and $\bvarpi_{m,q^a}$ with $\bomega_{\frac 12(a-m), \frac12(a+m)}$  we see that the irreducible objects in $\mathscr F_{n}$ are also indexed by elements of $\cal I_n^+$.\\\\
 Given $\bomega\in\cal  I^+_n$, the Weyl module $W(\bomega)$ is a universal finite-dimensional  cyclic $\widehat\bu_n$--module generated by an $\ell$--highest weight vector $v_{\bomegas}$; this means that, for each $1\leq i\leq n$ and $k\in \mathbb Z$ we have $x_{i,k}^+v_\bomegas=0$ and  $\phi_{i,k}^\pm$ acts on $v_\bomegas$ by a scalar determined by $\bomega$. 
Any  quotient of $W(\bomega)$ is called an $\ell$--highest weight module with $\ell$--highest weight $\bomega$ and it has  a unique irreducible quotient which is isomorphic to $V(\bomega)$. 
\\\\ 
   It was proved in  \cite{FR99} that  an object $V$  of $\mathscr F_{n}$ is  the  direct sum of generalized eigenspaces for the $\widehat{\bu}_n^0$--action. The eigenvalues are indexed by elements of $\cal I_n$ and we have, $$V=\bigoplus_{\bomegas\in\cal I_n} V_\bomegas,\ \ \wt_\ell V=\{\bomega\in\cal I_n: V_\bomegas\ne 0\},\ \ \wt_\ell ^\pm V=\wt_\ell V\cap (\cal I^+_n)^{\pm 1}.$$
     Moreover, if  $V'$ is another  object of $\mathscr F_{n}$ then  \begin{gather}
   \label{qinv}
   [V]=[V']\implies 
   \wt_\ell V=\wt_\ell V',\ \ \dim V_{\bomegas}= \dim V'_{\bomegas},\ \ \bomega\in\cal I_n,\\
\label{lwtprod}\wt_\ell (V\otimes V')=\wt_\ell V\wt_\ell V',\ \ \dim (V\otimes V')_\bomegas=\dim  (V'\otimes V)_\bomegas,\ \ \bomega\in\cal I_n.
   \end{gather}

\subsection{$\Sigma_r$-orbits, connected pairs, the maps $\tau_{m,\ell}$ and closed subsets  of $\mathbb I_n^r$}

\subsubsection{} The group $\Sigma_r$ acts on $\mathbb I_n^r$ in a natural way; for $\sigma\in\Sigma_r$ and $\bos\in\mathbb I_n^r$ set $$\sigma\bos=([i_{\sigma(1)}, j_{\sigma(1)}],\cdots,[i_{\sigma(r)}, j_{\sigma(r)}]).$$ The subsets $$\mathbb I_{n,\pm}^r=\{([i_1,j_1],\cdots, [i_r, j_r])\in\mathbb I_n^r: j_s-j_{s+1}\in\mathbb Z_\pm,\ \ 1\le s\le r-1\}$$ contain a set of orbit representatives for this action. Given $\bos\in\mathbb I_{n}^r$ set
$$n(\bos)+1=\max\{j_s: 1\le s\le r\}-\min\{i_p: 1\le p\le r\}.$$ Clearly $n(\bos)=n(w\bos)$ for all $w\in\Sigma_r$. Given $0\le p_1<p_2\le r$ set $$\bos(p_1,p_2)=([i_{p_1+1}, j_{p_1+1}],\cdots, [i_{p_2},j_{p_2}])\in\mathbb I_n^{p_2-p_1}.$$
Given $(\bos,\bos')\in \mathbb I_n^{r_1}\times \mathbb I_n^{r_2}$  we  write $\bos\vee\bos'$ for the  element of $\mathbb I_n^{r_1+r_2}$ obtained by concatenation. \\

\subsubsection{} Say that the pair  $([i_1, j_1], [i_2,j_2])\in\mathbb I_{n}^2$ is connected if either $i_2<i_1\le j_2<j_1$ and $0\le j_1-i_2\le n+1$ or $i_1<i_2\le j_1<j_2$ and $0\le j_2-i_1\le n+1$. Otherwise we say that they are not connected. 
\\\\ 
For $r\ge 1$, let
$\mathbb V_n^r$ be the vector space with basis $\mathbb I_n^r$. Define   $\tau:\mathbb V_n^2\to\mathbb V_n^2$  be the linear map defined on the basis  by   
\begin{gather*}\tau([i_1,j_1],[i_2,j_2])= ([i_2, j_1], [i_1,j_2]),\ \ {\rm if}\ \  ([i_1,j_1],[i_2,j_2] )\in\mathbb I_{n}^2\ \ {\rm is \ connected},\\
 \tau([i_1,j_1],[i_2,j_2])=0\ \ {\rm otherwise}.\end{gather*}
 Although $\tau$ obviously depends on $n$, we have for ease of notation omitted this dependence.\\\\
For $1\le m<\ell\le r$ let $\tau_{m,\ell}:\mathbb V_n^r\to\mathbb V_n^r$ be defined   by extending linearly the following assignment: if $([i_m,j_m],[i_\ell,j_\ell])$ is not connected we set $\tau_{m,\ell}\bos=0$ and, otherwise, we set 
$$\tau_{m,\ell}\bos = \bos(0,m-1)\vee ([i_\ell,j_m])\vee \bos(m,\ell-1)\vee ([i_m,j_\ell])\vee \bos(\ell,r).$$

\subsubsection{} Say that a subset $\mathbb J$ of $\mathbb I_{n}^r$ is closed if given  $\bos\in\mathbb J,\ 1\le m<\ell\le r$ we have
\begin{itemize}
\item $\tau_{m,\ell}\bos\in\mathbb J\sqcup\{0\}$,
\item and if $j_m=j_\ell$ then $\sigma_{m,\ell}\bos\in \mathbb J$.
\end{itemize} 
If $\mathbb J$ is not closed then we denote by $\bar{\mathbb J}[n]$ the smallest closed subset of $\mathbb I_{n}^r$ containing $\mathbb J$.  Clearly $\bar{\mathbb J}[n]\subset\bar{\mathbb{J}}[N]$ if $N\ge n$.
\\\\
If $\mathbb J=\{\bos\}$ then for ease of notation we denote its closure by $\bar\bos[n]$; if $\bar\bos[n]\subset \{w\bos: w\in\Sigma_r\}$ then we say that $\bos$ is closed.\\

\noindent{\bf Examples.}
\begin{itemize}
    \item
       The element  $\bos= ([0,6], [2,7], [1,8])$ is in $\mathbb I_{n,+}^3$ for any $n\ge 6$ and 
       \begin{gather*}\bar\bos[6]=\{\bos, ([2,6], [0,7],[1,8])\},\\ \bar\bos[n]=\bar\bos[6]\sqcup\{([1,6],[2,7],[0,8]), ([2,6], [1,7],[0,8])\}, \ \ n\ge 7.\end{gather*}
       \item Suppose that $\bos=([0,2], [1,2])$; then $\bar\bos[n]=\{\bos, ([1,2], [0,2])\}$, for $n\geq 2$. 
       \end{itemize}
       
\subsubsection{}  The proof of the next proposition can be found in Section \ref{comb}.
\begin{prop}\label{closedelement} Suppose that $\bos\in\mathbb I_{n,+}^r$ with $n\ge n(\bos)$. Then there exists $\bos_0\in\bar\bos[n]$ which is closed and the set of all closed elements in $\bar\bos[n]$ is contained in $\Sigma_r\bos_0$.
\end{prop}
\vskip12pt
\noindent{\bf Example.} We give an example to show that  the proposition is false if $n<n(s)$. Suppose that $n=1$ and take $\bos=([2,3], [1,2], [0,1])\in\mathbb I_{1,+}^3$ and note that $n(\bos)=2$. Then the elements $([1,3],[2,2], [0,1])$ and $([1,3], [1,1], [0,2])$ of $\bar\bos[1]$ are both closed and in distinct $\Sigma_3$ orbits.

 \subsection{Dominant $\ell$--weights}
 We state our first main result.
 \begin{thm} \label{weylmaps}
Let  $\bos\in\mathbb I_{n}^r$. 
\begin{enumerit}
\item[(i)] Given $\bos'\in\bar\bos[n]$, there exists an injective homomorphism $W(\bomega_{\bos'})\hookrightarrow W(\bomega_{\bos})$ of $\widehat\bu_n$--modules. 
\item[(ii)] We have 
$$\wt_\ell^+ W(\bomega_{\bos})=\{\bomega_{\bos'}: \bos'\in\bar\bos[n]\}.$$ 
In particular $W(\bomega_\bos)$ is irreducible if and only if $\bos$ is closed.
\end{enumerit}
 \end{thm}
 
 The following is immediate.
\begin{cor}
For $\bomega,\bpi\in\cal I_n^+$ we have
$$\Hom_{\widehat\bu_n}(W(\bpi), W(\bomega))\ne 0 \iff \bpi\in\wt^+_\ell W(\bomega).$$\end{cor}
\begin{rem} In Section \ref{tensor} we use the theorem to define new tensor subcategories of $\mathscr F_n$ which generalize the category $\mathscr C_\ell$ introduced in \cite{HL10}.
    
\end{rem}
\subsection{Maps between Weyl Modules and the Socle of a Weyl module}
Our second result is the following.
\begin{thm}\label{socleirr} Let $\bos\in\mathbb I_n^r$.\begin{enumerit}
    \item[(i)]
 We have $$\dim\Hom_{\widehat\bu_n}(W(\bomega_{\bos'}), W(\bomega_{\bos}))\le 1$$ with equality holding if and only if $\bos'\in\bar\bos[n]$.
\item[(ii)] Let $\bos_1,\cdots,\bos_k$ be the distinct closed elements (up to permutation) of $\bar\bos[n]$. 
Then,
$${\rm soc} (W(\bomega_\bos))\cong W(\bomega_{\bos_1})\oplus\cdots\oplus W(\bomega_{\bos_k}).$$
In particular, if  $n> n(\bos)$ then $k=1$ and $W(\bomega_\bos)$ has simple socle.
\end{enumerit}

    \end{thm}
    
\subsection{Some Applications}

 We give an application to the mixed Weyl modules. Given $\bos=([i_1,j_1],\cdots, [i_r,j_r])\in\mathbb I_n^r$, set$$W(\bos)=V(\bomega_{i_1,j_1})\otimes\cdots\otimes V(\bomega_{i_r, j_r}).$$
Then, it can be shown that 
$$W(\bos)\cong W(\bomega_{\bos(0,r_1)})\otimes\cdots\otimes W(\bomega_{\bos(r_{k-1},r_k)})$$ for some $k\ge 1$ and $1=r_0\le r_1\le\cdots\le r_k\le r$. 
An immediate corollary to Theorem \ref{weylmaps}(i) is 
 the following. Suppose that for $1\le p\le r$ we have $\bos_p'\in\overline{\bos(r_p,r_{p+1})}[n]\cap\mathbb I_{n,+}^r$. Setting $\bos'=\bos_1'\vee\cdots\vee \bos_k'$ we have an inclusion $ W(\bos')\hookrightarrow W(\bos)$.
\\\\
In Section \ref{mixedpf} we show that there are other kinds of inclusions between mixed Weyl modules as well. In particular, we are able to identify a particular irreducible module in the socle and in the  head of a mixed Weyl module.
\begin{prop}\label{mixed} Given $\bos\in\mathbb I_n^r$  there exist elements $\bos^\pm$ such that 
we have the following maps of $\widehat\bu_n$--modules:
$$V(\bomega_{\bos^-})\hookrightarrow W(\bos^-)\hookrightarrow W(\bos),\ \ \ W(\bos)\twoheadrightarrow W(\bos^+)\twoheadrightarrow V(\bomega_{\bos^+}).$$
\end{prop}
\noindent Finally, we also give a  sufficient condition in Section \ref{extensions} for  $\Ext^1_{\widehat\bu_n}(V,W)=0$ when  $V,W$ are objects of $\mathscr F_n$.

\section{The maps $\tau_{m,\ell}$ and the Proof of Proposition \ref{closedelement}}\label{comb}
\subsection{} We begin with the following property
\begin{lem}\label{braid} Suppose that $\bos=([i_1,j_1],\cdots, [i_r, j_r])\in\mathbb I_n^r$ is such  that the following hold:
$$ i_{s+1}<i_s\le j_r\  {\rm and} \  j_s<j_{s+1}, \  {\rm for} \   1\le s\le r-1.$$ 
Then for all $1\le m<p<\ell\le r$ with $j_m-i_\ell\le n+1$ we have
$$\tau_{m,p}\tau_{p,\ell}\tau_{m,p}\bos=\tau_{m,\ell}\bos.$$
\end{lem}
\begin{pf}
    It clearly suffices to prove the lemma when $m=1$ and $\ell=r$. Since $i_s<i_1\le j_s<j_1$ for $1<s\le r$ and $j_1-i_r\le n+1$ we see that 
   \begin{gather*}\tau_{1,r}\bos=([i_r, j_1])\vee\bos(1,r-1)\vee ([i_1,j_r]), \ \ {\rm and}\end{gather*}
   \begin{eqnarray*}
\tau_{1,p} \tau_{p,r}  \tau_{1,p}\bos&=&\tau_{1,0}\tau_{p,r}([i_p, j_1])\vee\bos(1,p-1)\vee ([i_1,j_p])\vee \bos(p,r)\\ 
&=&\tau_{1,p} ([i_p, j_1])\vee\bos(1,p-1)\vee([i_r,j_p])\vee \bos(p,r-1)\vee ([i_1,j_r])\\ 
&=&([i_r,j_1])\vee\bos(1,p)\vee\bos(p,r-1)\vee ([i_1,j_r]),
\end{eqnarray*}
where the final equality follows since $i_r<i_p<j_p<j_1$ and the lemma is proved. \end{pf}
\begin{rem}
It is in fact not hard to show that the maps $\tau_p:=\tau_{p,p+1}$, $1\le p\le r-1$, satisfy the following relations: for $1\le s,p\le r-1$,$$\tau_s^2=0,\ \ \tau_p\tau_s=\tau_s\tau_p,\ \ |s-p|\ge 2,\ \ \tau_p\tau_{p+1}\tau_p=\tau_{p+1}\tau_p\tau_{p+1}.$$ 
In particular the assignments $T_s\to \tau_s$,  $1\le s\le r-1$  define  an action of the Nil--Hecke algebra on  ${\rm {NH}}_r$  on $\mathbb V_k$. Denoting by $\tau_w$ the element of ${\rm {NH}}_r$ corresponding to $w\in \Sigma_r$ we have \begin{gather}\label{expltau}
\tau_w([i_1,j_1],\cdots, [i_r,j_r])\in\{0,([i_{w(1)}, j_{1}],\cdots ,[i_{w(r)}, j_{r}])\}.
\end{gather} 
  \end{rem} 
  \subsection{}
 For $\bos=([i_1,j_1],\cdots, [i_r,j_r])\in\mathbb I_{n,+}^r$ let $\Sigma_r(\bos)$ be the subgroup of $\Sigma_r$ generated by the transpositions $\sigma_{p}$ with $j_p=j_{p+1}$. We give an explicit description of $\bar\bos[n]$.

\begin{prop} Suppose that $\bos=([i_1,j_1],\cdots, [i_r,j_r])\in\mathbb I_{n,+}^r$. 
   For $1\le m<\ell\le r$ and $w\in\Sigma_r(\bos)$ we have\begin{gather*} \tau_{m,\ell} w\bos\ne 0\implies w^{-1}\tau_{m,\ell}w\bos= \tau_{w^{-1}(m), w^{-1}(\ell)}\bos,\\ 
  w\tau_{m,\ell}\bos\ne 0\implies w\tau_{m,\ell}\bos=  \tau_{w(m), w(\ell)}w\bos.
   \end{gather*}\end{prop}
\begin{pf} Since $j_1\ge j_2\ge \cdots\ge j_r$ it follows that $\Sigma_r(\bos)$ is generated by the simple transpositions $\sigma_p$, $1\le p\le r-1$ and  $j_p=j_{p+1}$.
    Hence  is enough to prove the proposition when $w=\sigma_{p}$ with $j_p=j_{p+1}$.\\\\
    The proposition is clear if $\{p,p+1\}\cap\{m,\ell\}=\emptyset$. If $m=p$ and $\ell=p+1$ then $$\tau_{m,\ell}\bos=0=\tau_{m,\ell}\sigma_p\bos$$ and there is nothing to prove.  It remains to   consider the case when  $\ell\in\{p,p+1\}$ and $m<p$ or $m\in\{p,p+1\}$ and $\ell>p+1$; we prove that
\begin{equation}\label{tauconj}\sigma_{\ell-\epsilon}\tau_{m,\ell}\sigma_{\ell-\epsilon}\bos=\tau_{m,\ell+1-2\epsilon}\bos,\ \ \ 
\sigma_{m-\epsilon}\tau_{m,\ell}\sigma_{m-\epsilon}\bos= \tau_{m+1-2\epsilon,\ell}\bos, \ \ \epsilon\in \{0,1\}.\end{equation}
We prove the first equality; the proof of the second is identical.\\\\
If $([i_m,j_m], [i_{\ell-\epsilon}, j_{\ell-\epsilon}])$ is not connected the we have $\tau_{m,\ell}\sigma_{\ell-\epsilon}=\tau_{m,\ell+1-2\epsilon}\bos=0$ and \eqref{tauconj} holds. Otherwise, if $\epsilon=0$ then 
$$\sigma_{\ell}\tau_{m,\ell}\sigma_\ell \bos = \bos(0,m-1)\vee ([i_{\ell+1},j_m])\vee \bos(m,\ell)\vee ([i_m,j_{\ell+1}])\vee \bos(\ell+1,r) = \tau_{m,\ell+1}\bos,$$
and  if $\epsilon = 1$ then 
$$\sigma_{\ell-1}\tau_{m,\ell}\sigma_{\ell-1}\bos = \bos(0,m-1)\vee ([i_{\ell-1},j_m])\vee \bos(m,\ell-2)\vee ([i_m,j_{\ell-1}])\vee\bos(\ell-1,r) = \tau_{m,\ell-1}\bos,$$ as needed. \\\\
Noting that $w\tau_{m,\ell}\bos= w\tau_{m,\ell}w^{-1}(w\bos)$ we have that  $\tau_{m,\ell}w^{-1}(w\bos)\neq 0$ if  $w\tau_{m,\ell}\bos\neq 0$. Since $w\bos\in \mathbb I_{n,+}^r$ it follows from the fist part that 
$$w\tau_{m,\ell}\bos =w\tau_{m,\ell}w^{-1}(w\bos) = \tau_{w(m),w(\ell)}w\bos,$$
as desired. 
\end{pf}
\begin{cor}\label{closure} Suppose that $\bos=([i_1,j_1],\cdots, [i_r, j_r])\in\mathbb I_{n,+}^r$.
Then 
$$\bar\bos[n]=\{w\tau_{m_1,\ell_1}\cdots \tau_{m_p,\ell_p}\bos: w\in\Sigma_r(\bos),\ \ p\ge 0,\  1\le m_s<\ell_s\le r\}\setminus\{0\}.$$
In particular, $\bar\bos[n]$ is a finite set and hence  contains closed elements of $\mathbb I_{n,+}^r$.
    \end{cor}
    \begin{pf}
        It suffices by the minimality of the closure to check that the set on the right-hand side is closed. But this is clear from the proposition.
    \end{pf}

 \subsection{}

 \begin{prop} \label{snakelike} Let $\bos=([i_1,j_1],\cdots, [i_r, j_r])\in\mathbb I_{n,+
 }^r$ with $n\ge n(\bos)$. There exists   $\sigma\in \Sigma_r$  with $i_{\sigma(1)}\ge \cdots\ge i_{\sigma(r)}$ such that $$ \bos_1=([i_{\sigma(1)}, j_1],\cdots, [i_{\sigma(r)}, j_r])\in\mathbb I_{n,+}^r\  {\rm and}\ \ \bos\in\bar\bos_1[n].$$
 
    \end{prop}
    \begin{pf} The proof proceeds by induction on $r$ with induction obviously beginning at $r=1$. For the inductive step choose $\sigma\in\Sigma_{r-1}$ such that $i_{\sigma(1)}\ge\cdots\ge i_{\sigma(r-1)}$ with 
    $$\bos_1'=([i_{\sigma(1)}, j_1],\cdots, [i_{\sigma(r-1)}, j_{r-1}])\in\mathbb I_{n,+}^{r-1},\ \ \bos(0,r-1)\in\overline{\bos_1'}[n].$$  
   We have  $\bos'=\bos_1'\vee ([i_r,j_r])\in\mathbb I_{n,+}^r$ and $\bos\in\overline{\bos'}[n]$. Hence it suffices to prove the proposition for $\bos'$; in other words we may assume without loss of generality that $\bos$ satisfies $i_1\ge\cdots\ge i_{r-1}$.
     If
    $i_r\le i_s$ for all $s\ge 1$ the proposition follows by taking $\bos_1=\bos'$.
    Otherwise set
    $$\bos_2=\bos(0,r-2)\vee([i_r,j_{r-1}])\vee([i_{r-1}, j_r]).$$ 
     We have \begin{gather*}
       j_r=j_{r-1} \implies \bos=\sigma_{r-1,r}\bos_2,\\
       j_r\ne j_{r-1} \implies i_{r-1}<i_r\leq j_r <j_{r-1}\implies \bos=\tau_{r-1, r}\bos_2.
      \end{gather*} In both cases we have $\bos\in\bar\bos_2[n]$. The inductive hypothesis applies to $\bos_2(0,r-1)$ and hence we can choose $\sigma'\in\Sigma_{r-1}$ such that the proposition holds for $\bos_2(0,r-1)$. 
      Taking $$\bos_1=([i_{\sigma(1)}, j_1],\cdots, [i_{\sigma'(r-1)}, j_{r-1}])\vee ([i_{r-1},j_r])$$ we get $\bos_2\in\bar\bos_1[n]$ and 
    the proposition follows. 
    \end{pf}
    
\subsection{Proof of Proposition \ref{closedelement}} In view of Proposition \ref{snakelike} it suffices to prove Proposition \ref{closedelement} under the assumption that $\bos\in\mathbb I_n^r$ satisfies both $i_s\ge i_{s+1}$ and $j_s\ge j_{s+1}$ for $1\le s\le r-1$.\\\\
Define
   $\sigma_\bos\in\Sigma_r$ as follows  \begin{gather*}{\sigma_\bos}(r)=\min\{s\in\{1,\cdots, r\}: i_s\le j_r\}\\
  \sigma_\bos(p)=\min\{s\in\{1,\cdots, r\}\setminus\{\sigma(p+1),\cdots, \sigma(r)\}: i_s\le j_p\}.
  \end{gather*} We prove first that $$\bos_0:=([i_{\sigma_\bos(1)},j_1],\cdots , [ i_{\sigma_\bos(r)}, j_r])\in\bar\bos[n],\ \ \bar\bos_0[n]\subset\{w\bos_0:w\in\Sigma_r\}.$$ 
It is clear from the definition that $i_{\sigma_\bos(s)}\le j_s$ for all $1\le s\le r$ and so $\bos_0\in\mathbb I_{n,+}^r$.  We check that $\bos_0$ is closed. 
For this, it suffices to prove the assertion that for all $1\le p<s\le r$ the pair $([i_{\sigma_\bos(p)}, j_p],[i_{\sigma_\bos(s)}, j_s])$ is not connected. This is clear if $j_s=j_p$ or $i_{\sigma_\bos(p)}=i_{\sigma_\bos(s)}$. Otherwise,  we have $i_{\sigma_\bos(s)}\le j_s<j_p$. It follows from the definition that either $i_{\sigma_\bos(p)}<i_{\sigma_\bos(s)}$ or $j_s<i_{\sigma_\bos(p)}$; in both cases the assertion is clear.
\\\\
We prove that $\bos_0\in\bar\bos[n]$ and that any closed element in $\bar\bos[n]$ is in $\Sigma_r\bos_0$ simultaneously. We proceed  by induction on $r$ with induction beginning at $r=1$.\\\\
Recall from Corollary \ref{closure} that $\bar\bos[n]$ must contain a closed element say $\bos_0'=([i_{w(1)}, j_1],\cdots , [i_{w(r)}, j_r])$ for some $w\in\Sigma_r$. If $i_{w(r)}= i_{\sigma_\bos(r)}$ then the induction hypothesis applies to $\bos(0,r-1)$ and gives $$\bos_0(0,r-1)\in\overline{\bos(0,r-1)}[n],\ \ \bos'(0,r-1)=w\bos_0(0,r-1)$$ for some $ w\in\Sigma_{r-1}$. The inductive step follows since we now also have $\bos_0'=w\bos_0$.
\\\\
If $i_{w(r)}\ne i_{\sigma_\bos(r)}$, let $1\le s<r$ be such that $w(s)=\sigma_\bos(r)$. The  definition of $\sigma_\bos$ forces $i_{w(r)}<i_{\sigma_\bos(r)}=i_{w(s)}\le j_r\le j_s$. Since $\bos_0'$ is closed this forces $j_r=j_{s}$. Hence $$\sigma_{s,r}\bos_0'=\bos_0'(0,s-1)\vee ([i_{w(r)}, j_r])\vee\bos_0'(s,r-1)\vee([i_{\sigma_\bos(r)},j_r])$$ is also a closed element of $\bar\bos[n]$ and so it suffices to prove that $\sigma_{s,r}\bos_0'\in \Sigma_r\bos_0$. For this, we
notice that $$\{i_{\sigma_\bos(1)},\cdots, i_{\sigma_\bos(r-1)}\}=\{i_{w(1)},\cdots ,i_{w(s-1)}, i_{w(r)},i_{w(s+1)},\cdots i_{w(r-1)}\}.$$ Choose $\sigma\in\Sigma_{r-1}$ such that 
$i_{\sigma\sigma_\bos(1)}\ge \cdots\ge i_{\sigma\sigma_\bos(r-1)}$, and set $$\bos_1=([i_{\sigma\sigma_\bos(1)}, j_1],\cdots,[i_{\sigma \sigma_\bos(r-1)}, j_{r-1}]).$$ Proposition \ref{snakelike} gives $$\{\bos_0(0,r-1), \ \ \sigma_{s,r}\bos_0'(0,r-1)\}\in\bar\bos_1[n].$$
Since $\bos_0(0,r-1)$ and $(\sigma_{r,s}\bos_0')(0,r-1)$ are both closed elements it follows from the inductive hypothesis that $(\sigma_{r,s}\bos_0')(0,r-1)=w'\bos_0(0,r-1)$ for some $w'\in\Sigma_{r-1}$. Since $w'(r)=r$ it follows that $\sigma_{r,s}\bos_0'=w\bos_0$ and the inductive step is proved.

\section{Proof of Theorem \ref{weylmaps}} \label{weylmapspf}

We begin by recalling two results which are needed for the proof.
\subsection{}\label{overlapred} The following is well-known, see for instance \cite{MY12}. \begin{lem} \label{lrootdrop}
    Suppose that $\bos=([i_1,j_1], [i_2,j_2])\in \mathbb I_{n}^2$ is connected. Then $$W(\bomega_{i_2,j_1}\bomega_{i_1,j_2})\cong W(\bomega_{\tau_{1,2}\bos})\hookrightarrow W(\bomega_\bos).$$ Moreover,
    \begin{gather*}
\wt^+_\ell W(\bomega_\bos)=\{\bomega_\bos,\ \bomega_{\tau_{1,2}\bos}\},\ \ \bomega_{i_1,j_1}^{-1}\bomega_{i_1,j_2}\bomega_{i_2,j_1}\in\wt_\ell V(\bomega_{i_2,j_2}).
\end{gather*}
\end{lem}

\subsection{}
The following result was  established in \cite{Ch01} (see also \cite{VV02}).
\begin{prop} \label{weylpermute} Suppose that $\bos=([i_1,j_1],\cdots, [i_k, j_k])\in\mathbb  I_n^k$. Then $$W(\bomega_{\bos})\cong V(\bomega_{i_1, j_1})\otimes \cdots\otimes V(\bomega_{i_k, j_k})$$ provided that  for all $1\le p<s\le k$ with $([i_p,j_p],[i_s,j_s])$ connected  we have $i_p+j_{p}\ge i_s+j_s$. 
 If $([i_p,j_p],[i_s,j_s])$ are not connected  for all $1\le s,p\le k$ 
 then $$W(\bomega)\cong V(\bomega)\cong V(\bomega_{i_{\sigma(1)}, j_{\sigma(1)}})\otimes \cdots\otimes V(\bomega_{i_{\sigma(k)}, j_{\sigma(k)}}),\ \ \sigma\in\Sigma_k.$$ 
 \end{prop}
  \begin{cor} Suppose that $\bos=([i_1,j_1],\cdots, [i_k,j_k])\in\mathbb I_n^k$ is such that $i_1\ge i_2\ge \cdots\ge i_k$ (or $j_1\ge\cdots\ge \ j_k$). Then
 $$W(\bomega_{\bos})\cong W(\bomega_{i_1,j_1})\otimes\cdots\otimes W(\bomega_{i_k, j_k}).$$
 \end{cor}
 \begin{pf} If $p<s$ and $[i_p, j_p]$  and $[i_s, j_s]$ is connected then  $i_s<i_p\le j_s<j_p$ and so the hypothesis of the proposition is satisfied. The corollary follows.
 \end{pf}

\subsection{Proof of Theorem \ref{weylmaps}(i)}  Since $\bomega_{\bos}=\bomega_{\sigma\bos}$ for all $\sigma\in\Sigma_r$ we can assume without loss of generality that $\bos\in\mathbb I_{n,+}^r$
Using Corollary \ref{closure} and the fact that $\bomega_{w\bos'}=\bomega_{\bos'}$ for all $w\in\Sigma_r$ and $\bos'\in\mathbb I_n^r$, it suffices to prove the result when  $$\bos'=\tau_{\ell,m}\bos,\ \ 
1\le \ell<m\le r\ \ {\rm and}\ \ ([i_\ell,j_\ell], [i_m,j_m]) \ \ {\rm is \ connected. }$$
We prove this by induction on $r$ with induction beginning at $r=2$ by Lemma \ref{overlapred}.
By Corollary \ref{weylpermute} we have 
\begin{gather*}W(\bomega_{\bos'})\cong W(\bomega_{\bos(0,\ell-1)})\otimes W(\bomega_{\bos'(\ell-1,m)})\otimes W(\bomega_{\bos(m,r)}),\\
W(\bomega_{\bos})\cong W(\bomega_{\bos(0,\ell-1)})\otimes W(\bomega_{\bos(\ell-1,m)})\otimes W(\bomega_{\bos(m,r)}).
\end{gather*}
If $\ell>1$ or $m<r$ then the inductive hypothesis gives an inclusion $W(\bomega_{\bos'(\ell-1,m)})\hookrightarrow W(\bomega_{\bos(\ell-1,m)})$ and hence the inductive step follows. \\\\
Hence we must prove  the inductive step when $\ell=1$ and $m=r$. In particular, we have that $([i_1,j_1],[i_r,j_r])$ is connected. Define a subset $\{p_1,\cdots, p_s\}$ of $\{1,\cdots, r\}$ as follows. Set $p_1=1$ and for $2\le m\le r$ define $p_m\le r$ to be minimal such that $([i_{p_{m-1}}, j_{p_{m-1}}], [i_{p_m}, j_{p_m}])$ is connected and let $s$ be maximal such that $p_s\le r$. In particular we have $i_{p_m}<i_{p_{m-1}}\le j_{p_m}<j_{p_{m-1}}$, for $2\leq m\leq s$. Together with the fact that $i_r<i_1\le j_r<j_1$ we get
\begin{equation}\label{connsn}i_{p_s}<i_{p_{s-1}}<\cdots<i_{p_2}<i_1\le j_r\leq j_{p_s}<j_{p_{s-1}}<\cdots<j_1.\end{equation}
If $p_s<r$ then $([i_{p_s}, j_{p_s}], [i_p,j_p])$ is not connected for all $p_s<p\le r$.  Setting 
$\bos''=\bos(0,p_s-1)\vee\bos (p_s, r)$, by Proposition \ref{weylpermute} we have $$W(\bomega_{\bos})\cong W(\bomega_{\bos''})\otimes W(\bomega_{i_{p_s}, j_{p_s}}).$$ 
The inductive hypothesis implies that we have a map $$W(\bomega_{\tau_{1,r-1}\bos''})\otimes W(\bomega_{i_{p_s}, j_{p_s}})\hookrightarrow W(\bomega_{\bos''})\otimes W(\bomega_{i_{p_s}, j_{p_s}})\cong W(\bomega_{\bos})$$ and hence it suffices to prove that $$W(\bomega_{\tau_{1,r-1}\bos''})\otimes W(\bomega_{i_{p_s}, j_{p_s}})\cong W(\bomega_{\tau_{1,r}\bos}).$$
By Corollary \ref{weylpermute}  we have $$W(\bomega_{\tau_{1,r-1}\bos''})\cong V(\bomega_{i_r,j_1})\otimes W(\bomega_{\bos(1,p_s-1)})\otimes W(\bomega_{\bos(p_s, r-1)})\otimes V(\bomega_{i_1,j_r}).$$ Since $i_{p_s}<i_1\le j_r<j_1$ it follows that $([i_{p_s}, j_{p_s}], [i_1,j_r])$
 is also not connected and 
 so Proposition \ref{weylpermute} gives $$W(\bomega_{\bos(p_s, r-1)})\otimes V(\bomega_{i_1,j_r})\otimes W(\bomega_{i_{p_s}, j_{p_s}})\cong W(\bomega_{i_{p_s}, j_{p_s}})\otimes W(\bomega_{\bos(p_s, r-1)})\otimes V(\bomega_{i_1,j_r}). $$
 Hence we get
 $$W(\bomega_{\tau_{1,r-1}\bos'})\otimes W(\bomega_{i_{p_s}, j_{p_s}})\cong V(\bomega_{i_r,j_1})\otimes W(\bomega_{\bos(1,p_s-1)})\otimes W(\bomega_{i_{p_s}, j_{p_s}})\otimes W(\bomega_{\bos(p_s, r-1)})\otimes V(\bomega_{i_1,j_r}).$$
 A further application of Corollary \ref{weylpermute} shows that the right hand side is isomorphic to $W(\bomega_{\tau_{1,r}\bos})$ as needed.\\\\
 If  $p_s=r$ then \eqref{connsn} shows that   Lemma \ref{braid} applies to $\tilde\bos=([i_1, j_1], [i_{p_2}, j_{p_2}],\cdots, [i_r, j_r])$ and gives 
 $$\tau_{1,s}(\tilde\bos)=\tau_{1,s-1}\tau_{s-1,s}\tau_{1,s-1}(\tilde\bos),\ \ {\rm and\ equivalently}\ \ \tau_{1,r}(\bos)=\tau_{1,p_{s-1}}\tau_{p_{s-1},r}\tau_{1,p_{s-1}}(\bos).$$
 Since $1<p_{s-1}<r$ it follows from the first part of the discussion, namely that $W(\bomega_{\tau_{m,\ell}\bos})\hookrightarrow W(\bomega_\bos)$ if $\ell-m<r-1$,
that we have the  following inclusions
$$W(\bomega_{\tau_{1,p_{s-1}}\tau_{p_{s-1},r}\tau_{1,p_{s-1}}(\bos)})\hookrightarrow W(\bomega_{\tau_{p_{s-1},r}\tau_{1,p_{s-1}}(\bos)})\hookrightarrow W(\bomega_{\tau_{1,p_{s-1}}(\bos)})\hookrightarrow W(\bomega_{\bos})$$ and the proof of part (i) of the theorem  is complete.

\subsection{}\label{mya} The proof of part (ii) requires   a very special case of a result established in   \cite{MY12}.\\\\ 
Given $[i,j]\in\mathbb I_n$ let $\mathbb P_{i,j}$ be the set of all functions $g:[0,n+1]\to\mathbb Z$ satisfying the following conditions:
\begin{gather*}g(0)=2j,\ \ g(r+1)-g(r)\in\{-1,1\},\ \ 0\le r\le n,\ \ g(n+1)=n+1+2i.\end{gather*}
 For $g\in\mathbb P_{i,j}$ we have $g(r)-r\in2\mathbb Z$ and we set
\begin{gather*}
\boc_g^\pm =\left\{\left[\frac12(g(r)-r), \frac12(g(r)+r)\right]: 1\le r \le n,\   g(r-1)=g(r)\pm 1=g(r+1)\right\},\\
\boc^\pm_{i,j}=\bigcup_{g\in \mathbb P_{i,j}}\boc_g^\pm,\\
\bomega(g)=\prod_{[m,\ell]\in\boc_g^+}\bomega_{m,\ell}\prod_{[m,\ell]\in \boc_g^-}\bomega_{m,\ell}^{-1}\in\cal I_n.
\end{gather*}  

It is straightforward to check that for $g\in \mathbb P_{i,j}$ we have 
\begin{gather}\label{uclc}
 [m,\ell]\in\boc_g^-\implies m+\ell>i+j.
 \end{gather}
 The next result is a particular case of \cite[Theorem 6.1]{MY12}.
\begin{prop}\label{mysnake}
\begin{enumerit}
    \item[(i)] Suppose that  $[i_1,j_1]\in \mathbb I_n$.
\begin{gather*}\label{dommon} \wt_\ell V(\bomega_{i_1,j_1})=\{\bomega(g): g\in\mathbb P_{i_1,j_1}\}, \ \ \dim V(\bomega_{i_1,j_1})_\bomegas \leq 1, \ \ \bomega\in\mathcal I_n^+,\\ \bomega(g)\in \wt_\ell^+ V(\bomega_{i_1,j_1})\iff \bomega(g)=\bomega_{i_1,j_1}\iff \boc_g^- = \emptyset.\end{gather*}
Moreover, $\bomega(g)\neq \bomega_{i_1,j_1}$ if and only if there  exists $[i,j]\in \boc_g^-$  such  that 
$$i+j > \max\{i'+j': [i',j']\in \boc_{g}^+\}\geq i_1+j_1.$$
\item[(ii)]  Suppose that $\bos = ([i_1,j_1],[i_2,j_2])\in \mathbb I_{n,+}^2$ is a connected pair and   $$\mathbb P_\bos=\{(g_1,g_2)\in \mathbb P_{i_1,j_1}\times \mathbb P_{i_2,j_2}: 
   g_{1}(m)>g_{2}(m),\ \ {\rm for\ all}\ \  0\leq m\leq n+1\}.$$ Then,
    $$\wt_\ell V(\bomega_\bos)=\{\bomega(g_1)\bomega(g_2): (g_1,g_2)\in\mathbb P_{\bos}\},\ \ \wt_\ell^+ V(\bomega_\bos)=\{\bomega_\bos\}.$$
    \end{enumerit}\hfill\qedsymbol    \end{prop}

\subsection{} The following was proved in \cite[Lemma 2.10]{BC25}.
\begin{lem}\label{ellwttau}
   Suppose that $\bos=([i_1,j_1], [i_2,j_2])\in\mathbb I_{n}^2$ with $i_1+j_1>i_2+j_2$. Then  $$[i_1,j_1]\in\boc_{i_2,j_2}^-\iff ([i_1,j_1], [i_2,j_2])\ \ {\rm is \ connected}.$$ \qed
\end{lem}

\subsection{Proof of Theorem \ref{weylmaps}(ii)} We prove part (ii), namely, $$\bomega\in\wt_\ell^+W(\bomega_{\bos})\iff  \bomega=\bomega_{\bos'},\ \ \bos'\in\bar\bos[n].$$ 
Part (i) of the theorem proves the reverse direction. We prove the forward direction
by induction on $r$.  The the case $r=2$ follows from Lemma \ref{lrootdrop}.  \\\\ 
Given $\bomega\in\wt^+_\ell W(\bomega_{\bos})$, we  use  equation \eqref{lwtprod}  and Proposition \ref{mysnake}(i)  to  write $$\bomega=\bomega(g_1)\cdots\bomega(g_r),\ \ g_s\in\mathbb P_{i_{s},j_{s}},\  \ 1\le s\le r. $$ 
Let $1\le p\le r$ be such that $i_p+j_p\ge\max\{i_s+j_s: 1\le s\le r\}$. We claim that $$\bomega(g_p)=\bomega_{i_p,j_p}\ \ {\rm and \ so}\ \ \bomega\bomega_{i_p,j_p}^{-1}\in\wt_\ell W(\bomega_\bos\bomega_{i_p,j_p}^{-1}).$$
Suppose for a contradiction that
$\bomega(g_p)\ne \bomega_{i_p,j_p}$. By Proposition \ref{mysnake}(i)  there exists $$[i,j]\in\boc_{g_p}^-\ \ {\rm with} \ \ i+j>i_{p}+j_{p}\ge i_{s}+j_{s},\ \ 1\le s\le r.$$  
On the other hand 
since $\bomega\in\cal I_n^+$ there exists $m\in[1,r]$ such that $[i,j]\in \boc_{g_m}^+$ and hence $\bomega(g_m)\neq \bomega_{i_m,j_m}$. Another application of Proposition \ref{mysnake}(i) implies that there  exists $[i',j']\in\boc_{g_m}^-$ with $i'+j'>i+j$. The argument iterates and contradicts $\bomega\in\cal I_n^+$ and the claim is proved.\\\\
The inductive hypothesis applies to $W(\bomega_\bos\bomega_{i_p,j_p}^{-1})$ and note that $\bomega_\bos\bomega_{i_p,j_p}^{-1}=\bomega_{\bos(0,p-1)\vee\bos(p,r)}$. Hence  if $\bomega\bomega_{i_p,j_p}^{-1}\in\cal I_n^+$ we get
$$\bomega\bomega_{i_p,j_p}^{-1}=\bomega_{\bos''},\ \ \bos''\in \overline{(\bos(0,p-1)\vee\bos(p,r))}[n].$$  
Since $\bos''(0,p-1)\vee ([i_p,j_p])\vee \bos''(p-1,r-1)\in \bar\bos[n]$ the inductive step is proved in this case.\\\\
If $\bomega\bomega_{i_p,j_p}^{-1}\notin\cal I_n^+$ then there exists $g_s$ with $s\ne p$ such that $[i_p,j_p]\in\boc_{g_s}^-$. It follows from \eqref{uclc} that $i_p+j_p>i_s+j_s$. Hence Lemma \ref{ellwttau}
gives that $([i_p,j_p], [i_s,j_s])$ is connected and in fact that $i_s<i_p\le j_s<j_p$. Since $\bos\in\mathbb I_{n,+}^r$ it follows that $p<s$ and $(g_p,g_s)\notin \mathbb P_{([i_p,j_p],[i_s,j_s])}$, since
$g_p(j_p-i_p)= i_p+j_p = g_s(j_p-i_p)$.
Moreover, we have $\bomega(g_p)=\bomega_{i_p,j_p}$ and $\dim V(\bomega_{i_s,j_s})_{\bomegas(g_s)}=1$, by Proposition \ref{mysnake}(i), and hence part (ii) implies that $\bomega(g_p)\bomega(g_s)\notin\wt_\ell V(\bomega_{i_p,j_p}\bomega_{i_s,j_s})$.  On the other hand, it follows from \eqref{lwtprod} that $$\bomega(g_p)\bomega(g_s)\in\wt_\ell (V(\bomega_{i_p,j_s})\otimes V(\bomega_{i_s, j_p}))=\wt_\ell W(\bomega_{i_s,j_p}\bomega_{i_p,j_s}).$$ Since 
\begin{gather*}\bomega\bomega(g_p)^{-1}\bomega(g_s)^{-1}\in\wt_\ell W(\bomega_\bos\bomega_{i_p,j_p}^{-1}\bomega_{i_s,j_s}^{-1}),\\ \wt_\ell (W(\bomega_\bos\bomega_{i_p,j_p}^{-1}\bomega_{i_s,j_s}^{-1})\otimes W(\bomega_{i_s,j_p}\bomega_{i_p,j_s}))= \wt_\ell W(\bomega_{\tau_{p,s}\bos}),\end{gather*}  we now have that 
$\bomega$ is an element of $\wt_\ell ^+ W(\bomega_{\tau_{p,s}\bos})$. Since $\overline{\tau_{p,s}\bos}[n]\subset\bar\bos[n]$, an iteration of this argument gives the result and the proof of part (ii) is complete.
Proposition \ref{weylpermute} was proved in greater generality in \cite{C01} and asserts that with this notation and assumptions we have $$W(\bvarpi)\cong W(\bvarpi_1)\otimes\cdots\otimes W(\bvarpi_k).$$

\section{Proof of Theorem \ref{socleirr}}\label{socle}
We begin with a discussion of the $\ell$--roots and some results on the $\ell$--weights of $V(\bomega_{i,j})$. The proof of the theorem starts in Section \ref{soclehoM}. \subsection{$\ell$--roots and a partial order on $\cal I_n^+$} \label{lrootdef}For $[i,j]\in\mathbb I_n$ with $0<j-i<n+1$
   set $$\balpha_{i,j}=\bomega_{i,j}\bomega_{i+1,j+1}(\bomega_{i+1,j}\bomega_{i,j+1})^{-1}.$$ Let $\cal Q_n^+$ be the submonoid (with unit) of $\cal I_n$ generated by the elements $\{\balpha_{i,j}: 0<j-i<n+1\}.$ It is well-known that $\cal Q_n^+$ is free on these generators and that if $\gamma\in\cal Q_n^+\setminus\{\bold 1\}$ then $
\gamma^{-1}\notin\cal I_n^+$. It is also clear that if  $m>n$ then we have a canonical inclusion $\cal Q^+_n\hookrightarrow \cal Q_m^+$.
Define a partial order $\preccurlyeq$ on $\cal I_n^+$ by $\bomega'\preccurlyeq\bomega $ iff $\bomega'=\bomega\balpha^{-1}$ for some $\balpha\in\cal Q_n^+$.
\\\\
The  proof of the next lemma is elementary and can be found in \cite{BC25}.
   \begin{lem}\label{gammaroot}
   Suppose that $([i_1,j_1], [i_2,j_2])\in\mathbb I_{n,+}^2$ is connected. Then$$\bomega_{i_1,j_1}\bomega_{i_2,j_2}=\bomega_{i_1,j_2}\bomega_{i_2, j_1}\prod_{i=i_2}^{i_1-1}\prod_{j=j_2}^{j_1-1}\balpha_{i,j}. $$ \hfill\qedsymbol
   \end{lem}
The following is well-known (see \cite{Ch01} for instance). Suppose that $[i,j]\in\mathbb I_n$; then $$\bomega_{j,n+1+i}^{-1}\in\wt_\ell V(\bomega_{i,j})\ \ {\rm and} \ \ \bomega\in\wt_\ell V(\bomega_{i,j})\implies \bomega_{j,n+1+i}^{-1}\preccurlyeq\bomega\preccurlyeq \bomega_{i,j}.$$
Since $([j,n+1+i], [i,j])\in\mathbb I_{n,+}^2$ is connected it follows from Lemma \ref{gammaroot} that \begin{gather}\label{sqfree}\gamma\in\cal Q_n^+,\ \ \bomega=\bomega_{i,j}\gamma^{-1}\in\wt_\ell V(\bomega_{i,j})\implies \gamma^{-1}\prod_{s=i}^{j-1}\prod_{m=j}^{n+i}\balpha_{s,m}\in\cal Q_n^+\\\notag \implies \gamma\balpha_{p,m}^{-2}\notin\cal Q_n^+\ \ {\rm for \ all}\  [p,m]\in\mathbb I_n.
\end{gather}
Moreover it also follows that if 
\begin{gather}\label{sqfree2}\bomega\in\wt_\ell V(\bomega_{i,j}), \ \ {\rm with} \ \ \bomega\preccurlyeq\bomega_{i,j}\eta^{-1},\ \ \eta=\prod_{s=i}^{j-1}\prod_{m=j}^{j_1},\ \ j-i\le j_1-i\le n,
\end{gather} then $\bomega=\bomega_{i,j}\eta^{-1}\gamma^{-1}$ where $\gamma$ is in the submonoid of $\cal Q_n^+$ generated by the elements $\balpha_{m,\ell}$ with $\ell\ge j+1$.

\subsection{The case $n>n(\bos)$}\label{soclehoM} We  establish Theorem \ref{socleirr} in the case when $n>n(\bos)$.
\subsubsection{} The first step of the proof is the following.
\begin{prop}\label{dim1}
    Suppose that $\bos=([i_1,j_1],\cdots, [i_r, j_r])\in\mathbb I_{n,+}^r$ is such that  $i_s\ge i_{s+1}$ for all $1\le s\le r-1$ and $n>n(\bos)$.
    Setting  $$\bomega=(\bomega_{i_1,j_1+1}\bomega_{j_1,j_1+1}^{-1})(\bomega_{i_2,j_1+1}\bomega_{j_2,j_1+1}^{-1})\cdots(\bomega_{i_r,j_1+1}\bomega_{j_r,j_1+1}^{-1}),$$
    we have  $$\dim W(\bomega_{\bos'})_{\bomegas}=1\ \ {\rm  for\ all} \ \ \bos'\in\bar\bos[n].$$
\end{prop}
\begin{pf} By Corollary \ref{closure} and \eqref{expltau} we can write $\bos'=([i_{\sigma(1)}, j_1],\cdots, [i_{\sigma(r)}, j_r])$ for some  $\sigma\in\Sigma_r$. Then we have
$$\bomega_{j_s, j_1+1}^{-1}\bomega_{i_{\sigma(s)}, j_1+1}\in\wt_\ell V(\bomega_{i_{\sigma(s)}, j_s}),\ \  1\le s\le r.$$
If  $i_{\sigma(s)}=j_s$ the statement holds vacuously and otherwise we have that $([j_s,j_1+1],[i_{\sigma(s)}, j_s])\in\mathbb I_{n,+}^2$ is connected and hence the statement follows from Lemma  \eqref{lrootdrop}. Since we can also write 
$$\bomega =(\bomega_{i_{\sigma(1)}, j_1+1}\bomega_{j_1,j_1+1}^{-1})\cdots(\bomega_{i_{\sigma(r)}, j_1+1}\bomega_{j_r,j_1+1}^{-1}),$$
it follows from \eqref{lwtprod} that $\dim W(\bomega_{\bos'})_\bomegas\ge 1$, for all $\bos'\in \bar\bos[n]$.
By Theorem \ref{weylmaps}(i) we have an inclusion $W(\bomega_{\bos'})\hookrightarrow W(\bomega_\bos)$ and hence  the proposition follows if we  prove that $\dim W(\bomega_\bos)_\bomegas=1$. We proceed by induction on $r$ with induction beginning at $r=1$. Setting $$\bos_1=([j_1,j_1+1], \cdots, [j_r, j_1+1])\in\mathbb I_{n,+}^r,\  \bos_2=([i_1,j_1+1], \cdots, [i_r, j_1+1])\in\mathbb I_{n,+}^r,$$ and using Lemma \ref{gammaroot} we see that 
$$\bomega =\bomega_{\bos_2}\bomega_{\bos_1}^{-1} =\bomega_{\bos}\eta^{-1},\ \ \eta=\eta_1\cdots \eta_r,\ \ \eta_s=\prod_{i=i_s}^{j_s-1}\prod_{j=j_s}^{j_1}\balpha_{i,j},\ \ 1\le s\le r.$$
 For $i_r\le i<j\le j_1$ we have $$\eta=\prod_{i_r\le i<j\le j_1}\balpha_{i,j}^{p_{i,j}},\ \ p_{i,j}=\#\{1\le s\le r: \ i_s\le i<j_s\le j\le j_1\}.$$ 
 By \eqref{lwtprod} we can  write $$\bomega=\bomega_\bos\eta^{-1}=\bomega_1\bomega_2,\ \ \bomega_1\in \wt_\ell W(\bomega_{\bos(0,r-1)}),\ \ \bomega_2\in\wt_\ell W(\bomega_{i_r,j_r}).$$
Suppose that $i_r\le i<j_r\leq j\le j_1$. Then $$p_{i,j}= p_{i,j}'+1,\ \ p_{i,j}'=\#\{1\le s\le r-1: i_s\le i<j_s\le j\le j_1\},$$ and so it  follows from equation \eqref{sqfree} that $\bomega_2\preccurlyeq\bomega_{i_r,j_r}\eta_r^{-1}$. In other words $\bomega_2$ satisfies the conditions of \eqref{sqfree2} and so $\bomega=\bomega_{i,j}\eta_r^{-1}\gamma^{-1}$ where  $\gamma$ is in the submonoid generated by $\balpha_{\ell, m}$ with $m>j_1$. Since $\eta$ has no such terms this forces $\bomega=\bomega_{i,j}\eta_r^{-1}$. Hence we have proved that 
$$\bomega_2=\bomega_{i_r, j_r}\eta_r^{-1} \ \ {\rm  and\  so }\ \ \bomega_1=\bomega_{\bos(0,r-1)}(\eta_1\cdots\eta_{r-1})^{-1}\in\wt_\ell W(\bomega_{\bos(0,r-1)}).$$ It follows that 
$$W(\bomega_{\bos})_\bomegas\cong W(\bomega_{\bos(0,r-1)})_{\bomegas_1}\otimes W(\bomega_{i_r,j_r})_{\bomegas_{i_r,j_r}\eta_r^{-1}},$$
and hence the inductive hypothesis applies and gives the inductive step.
\end{pf}

\begin{cor}\label{soclehom} Let $\bos\in\mathbb I_{n,+}^r$ and let  $\bos_0$ be  the unique closed element in $\bar\bos[n]$.  Then $$\dim\Hom_{\widehat{\bu}_n}(W(\bomega_{\bos_0}), W(\bomega_\bos))=1$$ and any non-zero element of this space is injective.
\end{cor}
\begin{pf} By Theorem \ref{weylmaps}(i) we have that $\dim\Hom_{\widehat{\bu}_n}(W(\bomega_{\bos_0}), W(\bomega_\bos))\ge 1$ and also that there exists an injective map $W(\bomega_{\bos_0})\hookrightarrow W(\bomega_\bos)$.
 By Theorem \ref{weylmaps}(ii), we have that $W(\bomega_{\bos_0})$ is irreducible and so any non-zero map $\eta: W(\bomega_{\bos_0})\to W(\bomega_\bos)$
must be injective. In particular it is non-zero on $W(\bomega_{\bos_0})_{\bomegas}$ where $\bomega$ is as in the  proposition. Since this $\ell$--weight space generates $W(\bomega_{\bos_0})$ it follows   that $\eta$ is determined by its value on this one dimensional space and so $\dim\Hom_{\widehat{\bu}_n}(W(\bomega_{\bos_0}), W(\bomega_\bos))=1$.
   \end{pf}

\subsubsection{}\label{duals} We recall some well-known results on duals. Any object $V$ of $\mathscr F_n$  has a right and a  left dual  denoted by $V^*$ and ${}^*V$ respectively, and we have $\widehat{\bu}_n$--maps
 $$\mathbb C\hookrightarrow V^*\otimes V,\ \ {}^*V\otimes V\to \mathbb C\to 0.$$
 Moreover, if $\bos = ([i_1,j_1],\cdots, [i_r,j_r]) \in \mathbb I_n^r$ then $V(\bomega_\bos)^*\cong V(\bomega_{\bos^*})$ and ${}^*V(\bomega_\bos)\cong V(\bomega_{{}^*\bos})$,
 where
 $$\bos^* = ([j_1,n+1+i_1],\cdots, [j_r,n+1+i_r]), \ \  {}^*\bos = ([-n-1+j_1,i_1],\cdots, [-n-1+j_r, i_r]).$$
We shall freely use properties of duals, in particular, the isomorphisms
\begin{gather*}(W\otimes V)^*\cong V^*\otimes W^*,\ \  {}^*(W\otimes V)\cong {}^*V\otimes {}^*W,\\
 \Hom_{\widehat{\bu}_n}(V\otimes U, W)\cong \Hom_{\widehat {\bu}_n}( U, V^*\otimes W),\ \ \Hom_{\widehat{\bu}_n}(U\otimes V, W)\cong \Hom_{\widehat{\bu}_n}( U, W\otimes {}^*V).\end{gather*}

\subsubsection{} We prove Theorem \ref{socleirr}(i); by Theorem \ref{weylmaps}(i) it suffices to prove that $$\dim \Hom_{\widehat{\bu}_n}(W(\bomega_{\bos'}), W(\bomega_{\bos}))\le 1,\ \ \bos,\bos'\in\mathbb I_{n,+}^r,\ \ n>n(\bos).$$ 
Let $
\bos_1=([i_{\sigma(1)},j_1],\cdots, [i_{\sigma(r)}, j_r)]$ be such that $i_{\sigma(s)}\ge i_{\sigma(s+1)}$ for $1\le s\le r-1$. Then $\bos\in\bar\bos_1[n]$ by Proposition \ref{snakelike} and we have an inclusion $W(\bomega_{\bos})\hookrightarrow  W(\bomega_{\bos_1})$ by Theorem \ref{weylmaps}(i). Hence, 
$$\dim\Hom_{\widehat{\bu}_n}(W(\bomega_{\bos'}), W(\bomega_{\bos}))\le \dim \Hom_{\widehat{\bu}_n}(W(\bomega_{\bos'}), W(\bomega_{\bos_1})), $$ 
and it suffices to prove the assertion  for $\bos_1$; or equivalently when $\bos=([i_1, j_1],\cdots, [i_r,j_r])$ satisfies $i_s\ge i_{s+1}$ and $j_s\ge j_{s+1}$ for $1\le s\le r-1$ and we assume this is the case from now on.
\\\\
We proceed by induction on the length of a reduced expression for $\bomega_{\bos'}$.
Suppose that  $\bomega_{\bos'}=\bomega_{i,j}$ for some $[i,j]\in\mathbb I_n$ and that $0\ne \eta\in \Hom_{\widehat{\bu}_n}(W(\bomega_{i,j}), W(\bomega_{\bos}))$. Since $([i,j])$ is closed  Proposition \ref{closedelement} implies it must be the unique closed element of $\bar\bos[n]$ and hence Corollary \ref{soclehom} shows that induction begins.\\\\
Assume the result for all $\bos'$ for which the  length of a reduced expression for $\bomega_{\bos'}$ is $N-1$.    Suppose that $\bos'$ is such that the length of a reduced expression of $\bomega_{\bos'}$ is $N$ and that $\Hom_{\widehat{\bu}_n}(W(\bomega_{\bos'}), W(\bomega_{\bos}))\ne 0$. Then by Theorem \ref{weylmaps}(ii) and its corollary we have 
$$\bos'=([i_{p_1},j_1],\cdots, [i_{p_r}, j_r])\in\bar\bos[n],$$
where $\{p_1,\cdots, p_r\}$ is a permutation of $\{1,\cdots, r\}$. Moreover, we can also assume that $p_1<p_s$ if $j_s=j_1$.
Since $i_{p_1}\le i_1<j_1$ and $j_1-i_{p_1}\le n$ we have 
 that $\bomega_{i_{p_1}, j_1}\ne \bold 1$.
If  $\bomega_{i_\ell, j_m}\ne\bold 1$ also occurs in a reduced expression of $\bomega_{\bos'}$ then our choice give that either $i_{p_1}+j_1\ge i_\ell+j_m$ or 
$i_{p_1}\le i_\ell<j_m\le j_1$. In either case Proposition \ref{weylpermute} implies that $$(*)\ \ 
W(\bomega_{\bos'})\cong V(\bomega_{i_{p_1},j_1})\otimes W(\bomega_{\bos'}\bomega_{i_{p_1},j_1}^{-1})$$ and so $$(**)\ \ \Hom_{\widehat{\bu}_n}(W(\bomega_{\bos'}), W(\bomega_{\bos}))\cong\Hom_{\widehat{\bu}_n}(W(\bomega_{\bos'}\bomega_{i_{p_1},j_1}^{-1}), V(\bomega_{j_1,n+1+i_{p_1}})\otimes W(\bomega_{\bos})).$$
Since $i_1-i_{p_1}< j_1-i_{p_1}\le n$ we have $n+i_{p_1}+j_1\ge j_1+i_1\ge j_s+i_s$ for all $1\le s\le r$. Another application of Proposition \ref{weylpermute} now gives $$V(\bomega_{j_1,n+1+i_p})\otimes W(\bomega_{\bos})\cong W(\bomega_{j_1,n+1+i_p}\bomega_{\bos}).$$ Since the length of a reduced expression for $\bomega_{\bos'}\bomega_{i_p,j_1}^{-1}$ is $N-1$, the  inductive hypothesis applies to the pair $\bomega_{\bos'}\bomega_{i_p,j_1}^{-1}$ and $\bomega_{([j_1,n+1+i_p])\vee\bos}$.
It follows that the right hand side of $(**)$ is one-dimensional and hence the inductive step is proved. 
\subsubsection{} We prove Theorem \ref{socleirr}(ii).
\\\\
 Suppose that we have a non-zero map $V(\bpi)\to W(\bomega_\bos)$ and hence a non-zero composite map $\eta: W(\bpi)\to V(\bpi)\to W(\bomega_\bos)$. Then by Theorem \ref{weylmaps}(ii)  we have $\bpi=\bomega_{\bos'}$ for  some $\bos'\in\bar\bos[n]$. By  Theorem \ref{socleirr}(i) with $n>n(\bos)$ it follows that $\eta$ must be injective  and so $W(\bomega_{\bos'})\cong V(\bomega_{\bos'})$. But this means by Theorem \ref{weylmaps}(i) that $\bos'$ must be closed and so by Proposition \ref{closedelement} that $\bos'\in \Sigma_r\bos_0$. Hence $W(\bomega_{\bos_0})$ is the socle of $W(\bomega_\bos)$.

\subsection{The case $n\le n(\bos)$.} We now deduce the analogous results for all $\bos\in\mathbb I_n^r$ for general $n$. 
Recall that in this case the set $\bar\bos[n]$ could contain more than one closed element. \\\\
For  $N>n(\bos)$ and given $\bos\in\mathbb I_n^r$ we denote by $W_N(\bomega_\bos)$ the Weyl module for $\widehat\bu_N$. The algebra  $\widehat\bu_n$ maps to  the subalgebra $\tilde\bu_n$ of $\widehat\bu_N$ generated by elements $x_{i,k}^\pm$, $\mathbb\phi_{j,k}^\pm$ with $1\le i\le n$, $1\le j\le N$ and $k\in\mathbb Z$. 
We have a surjective map $\iota: W_n(\bomega_\bos)\to \tilde\bu_n w_{\bos}$ of $\widehat\bu_n$--modules where $w_\bos$ denotes the $\ell$--highest weight vector of $W_N(\bomega_{\bos})$. We prove first that this map must be an isomorphism.
\\\\
For this,  observe that we have isomorphisms of vector spaces 
\begin{eqnarray*}\tilde\bu_n w_\bos &=&\bigoplus_{\eta\in \cal Q_n^+} W_N(\bomega_\bos)_{\bomegas_\bos\eta^{-1}}\\ &=&\bigoplus _{\eta_1,\cdots,\eta_r\in\cal Q_n^+}\left(W(\bomega_{i_1,j_1})_{\bomegas_{i_1,j_1}\eta_1^{-1}}\otimes \cdots\otimes W(\bomega_{i_r,j_r})_{\bomegas_{i_r,j_r}\eta_r^{-1}}\right) \\ &\cong& W_n(\bomega_{i_1,j_1})\otimes\cdots\otimes W_n(\bomega_{i_r, j_r}).\end{eqnarray*} It follows that $\dim W_n(\bomega_\bos)=\dim\tilde\bu_nw_\bos$ and hence $\iota$ is an isomorphism.\\\\
Suppose now that we have a non-zero map $\eta: W_n(\bomega_{\bos'})\to W_n(\bomega_\bos)$ and hence a non-zero map
$$\tilde \eta: W_n(\bomega_{\bos'})\to W_n(\bomega_\bos)\cong\tilde\bu_n w_\bos\hookrightarrow W_N(\bomega_\bos).$$ Then it follows from \cite{FM01} that $\tilde\eta(w_{\bos'})$ is an $\ell$--weight vector with $\ell$--weight $\bomega_{\bos'}$. Moreover since $w_{\bos'}$ is $\ell$--highest weight it follows $x_{i,k}^+\tilde\eta(w_{\bos'})=0$ for all $1\le i\le n$ and $k\in\mathbb Z$. Since $\bomega_{\bos}\bomega_{\bos'}^{-1}$ is in the subgroup of $\cal Q_N^+$ generated by $\balpha_{i,j}$ with $0<j-i<n+1$ it is also true that $x_{i,k}^+\tilde\eta(w_{\bos'})=0$
for all $n< i\le N$ and $k\in\mathbb Z$. Hence we have a non-zero map $W_n(\bomega_{\bos'})\hookrightarrow W_N(\bomega_{\bos'})\hookrightarrow W_N(\bomega_{\bos}) $. Since Theorem \ref{socleirr} hold if $N>n(\bos)$ it now follows that 
it also holds for $n\le n(\bos)$.
\\\\
Next we prove that the socle of $W_n(\bomega_{\bos})$ is isomorphic to the direct sum of $W_n(\bomega_{\bos_0})$ (with multiplicity one) where $\bomega_{\bos_0}$ varies over all the closed elements of $\bar\bos[n]$. Suppose that we have a  map  of $\widehat{\bu}_n$--modules $V(\bpi)\to W_n(\bomega_{\bos})$. Composing with $W_n(\bpi)$ we get a non-zero map $W_n(\bpi)\to V(\bpi)\to  W_n(\bomega_\bos)$ 
which has to be injective. Hence $\bpi=\bomega_{\bos'}$ for some $\bos'\in\bar\bos[n]$ by Theorem \ref{weylmaps}(ii) and $W_n(\bomega_{\bos'})\cong V(\bomega_{\bos'})$. It follows that $\bos'$ is a closed element of $\bar\bos[n]$ and the proof is complete.

\section{Applications} \label{mixedpf}
In this section we give three applications of our results, including the proof of Proposition \ref{mixed}.
\subsection{Tensor Subcategories} \label{tensor} We begin with  application of Theorem \ref{weylmaps}. \\\\
Given $\bos=([i_1,j_1],\cdots, [i_r,j_r])\in\mathbb I_{n}^r$ let $\mathcal I_n^+[\bos]$ be the submonoid of $\cal I_n^+$ generated by elements $\bomega_{i_s,j_p}$ for $1\le s,p\le r$ with $0\le j_p-i_s\le n+1$. We have \begin{gather*}\bomega\in\cal I_n^+[\bos]\implies\bomega=\bomega_{\bos'},\ \ \bos'=([m_1,\ell_1],\cdots, [m_k,\ell_k]),\ \ k\ge 1,\\
{\rm and\  for}\ \ 1\le s\le k\ \ {\rm we \ have}\ \ m_s=i_p,\ \ \ell_s=j_{p'}\ \ {\rm for \ some}\ \ 1\le p,p'\le r.\end{gather*} 
It follows that \begin{equation}\label{ins}\bos''\in\bar{\bos}'[n]\implies \bomega_{\bos''}\in\cal I_n^+[\bos].\end{equation}
\begin{prop} Let $\mathscr F_n(\bos)$ be the full subcategory of $\mathscr F_n$ consisting of objects whose Jordan Holder constituents are of the form $[V(\bomega)]$ with $\bomega\in\cal I_n^+[\bos]$. Then $\mathscr F_n[\bos]$ is an abelian tensor subcategory.
\end{prop}
\begin{pf}
Suppose that $V(\bomega)$ occurs in the Jordan--Holder series of $V(\bomega_{\bos_1})\otimes V(\bomega_{\bos_2})$.  Then $$\bomega\in\wt_\ell^+(W(\bomega_{\bos_1})\otimes W(\bomega_{\bos_2}))\iff\bomega\in\wt_\ell^+W(\bomega_{\bos_1}\bomega_{\bos_2}).$$
 Writing $\bomega=\bomega_{\bos'}$ for some $\bos'\in\mathbb I_n^k$, $k\ge 1$  and noting that $\bomega_{\bos_1}\bomega_{\bos_2}=\bomega_{\bos_1\vee\bos_2}$, it follows from  Theorem \ref{weylmaps}(ii)  that $$\bos'\in\overline{\bos_1\vee\bos_2}[n].$$ It is immediate from Corollary \ref{closure} and \eqref{expltau} that $\bomega_{\bos'}\in\cal I_n^+[\bos]$  and the proof of the proposition is complete.
\end{pf}
\begin{rem}
   The categories $\mathscr F(\bos)$  generalize the categories $\mathscr C_\ell$ of \cite{HL10}. In fact our   results  give a different proof of the fact that the various subcategories considered in  \cite{HL10} and \cite{HL13a} are tensor subcategories.
\end{rem}

 Notice that if  $\bos\in\mathbb I_{n,+}^r$ then Proposition \ref{mixed} is a special case of Theorem \ref{socleirr} and Corollary \ref{weylpermute} where we take $\bos^-$ to be a closed element of $\bar\bos[n]$ and $\bos^+=\bos$.

\subsection{}  We establish  Proposition \ref{mixed} in the case when $\bos\in\mathbb I_{n,-}^r$.
\begin{prop}\label{antiweyl}
  Suppose that $\bos\in\mathbb I_n^r$ is such that for all $1\le p< s\le r$ with $([i_p,j_p], [i_s,j_s])$  connected we have $i_s+j_s\ge i_p+j_p$. Let $\bos_0$ be a closed element of $\bar\bos[n]$. We have the following maps: $$V(\bomega_{\bos})\hookrightarrow W(\bos)\twoheadrightarrow V(\bomega_{\bos_0}).$$
\end{prop}
\begin{pf}
We start by noting that $\overline{^*\bos}[n]=\{^*\bos':\bos'\in\bar\bos[n]\}$. In fact, it suffices to show that $([i,j],[i',j'])$ is connected if and only if ${}^*([i,j],[i',j'])$ is connected. If $j\leq j'$, it follows by noting that 
$$ i<i'\leq j<j' \ {\rm and}\  j'-i\leq n+1\iff -n-1+j<-n-1+j'\leq i<i' \ \ {\rm and} \ i'\leq j.$$
The case $j'<j$ is clearly analogous. \\

Let $\bos_0\in \bar\bos[n]$ be closed. Then ${}^*\bos_0\in \overline{{}^*\bos}[n]$ is closed and then Proposition \ref{weylpermute} and Theorem \ref{socleirr} give $$W(\bomega_{^*\bos_0})\hookrightarrow W(\bomega_{^*\bos})\twoheadrightarrow V(\bomega_{^*\bos}).$$
Taking duals and noting that 
$W(\bomega_{^*\bos})^*\cong W(\bos)$  gives the result.
\end{pf}

\subsection{The maps $\iota^\pm$} Define $\iota^\pm:\mathbb I_n^2\to\mathbb I_{n,\pm}^2$ as follows:
\begin{gather*}\iota^+([i_1,j_1],[i_2, j_2])=\begin{cases} ([i_1,j_1],[i_2,j_2]),\ \ j_1\ge j_2,\\([i_2,j_2],[i_1,j_1]),\ \ j_2>j_1\  {\rm and}\  ([i_1,j_1], [i_2,j_2])\ {\rm is\ not\ connected},\\([i_1,j_2], [i_2,j_1]),\ \  j_2>j_1\  {\rm and} \ ([i_1,j_1], [i_2,j_2])\  {\rm is\  connected},\end{cases}\\
\iota^-([i_1,j_1],[i_2, j_2])=\begin{cases} ([i_1,j_1],[i_2,j_2]),\ \ j_2\ge j_1,\\
([i_2,j_2],[i_1,j_1]),\ \ j_1>j_2\  {\rm and} \ ([i_1,j_1], [i_2,j_2])\ {\rm is\ not\ connected},\\
([i_1,j_2], [i_2,j_1]),\ \  j_1>j_2\  {\rm and}\ ([i_1,j_1], [i_2,j_2])\ {\rm is\ connected}.\end{cases}
\end{gather*}
Notice that $\sigma\iota^+\sigma=\iota^-$ where $\sigma:\mathbb I_n^2\to\mathbb I_n^2$ is the permutation.
\\\\
\noindent {\bf Example.} Suppose that $n\geq 7$. Then
\begin{gather*}
\iota^+([2,5],[3,9])= ([2,9],[3,5]), \ \ \iota^+([3,9],[2,5])= ([3,9],[2,5]),\ \ \iota^+([3,5],[2,9]) = ([2,9],[3,5]).
\end{gather*}

\subsection{} Given $\bos\in\mathbb I_n^r$ and $1\le p\le r-1$ set $$\iota_p^\pm\bos=\bos(0,p-1)\vee\iota^\pm(\bos(p-1,p+1))\vee\bos(p+1,r).$$ 

\begin{prop}
    For $\bos\in\mathbb I_n^r$ and $1\le p\le r-1$  we have $$W(\iota_p^-\bos)\hookrightarrow W(\bos)\twoheadrightarrow W(\iota^+_p\bos).$$
\end{prop}
\begin{pf}
    If $ j_{p+1}\ge j_p$ or if $([i_p,j_p], [i_{p+1}, j_{p+1}])$ is not connected then by Proposition \ref{weylpermute} we have $W(\iota_p^-\bos)\cong W(\bos)$. If $j_p>j_{p+1}$ and $([i_p,j_p], [i_{p+1}, j_{p+1}])$ is connected then by Lemma \ref{lrootdrop} we have a map $$W(\iota^-\bos(p-1,p+1))\hookrightarrow W(\bos(p-1,p+1))\ \  {\rm and \ so}$$
    $$W(\iota_p^-\bos)\cong W(\bos(0,p-1)\otimes W(\iota^-\bos(p-1,p+1))\otimes W(\bos(p+1,r))\hookrightarrow W(\bos).$$
    Similarly  if $j_p\ge j_{p+1}$ or if $([i_p,j_p], [i_{p+1}, j_{p+1}])$ is not connected we have $W(\iota^+\bos)\cong W(\bos)$. 
    Otherwise by Lemma \ref{lrootdrop}  we have $$W(\bos(p-1, p+1))\twoheadrightarrow W(\iota^+\bos(p-1,p+1))$$ and hence 
    $$W(\bos)\cong W(\bos(0,p-1))\otimes W(\bos(p-1,p+1))\otimes W(\bos(p+1,r))\twoheadrightarrow W(\iota_p^+\bos). $$
\end{pf}
The following is immediate.
\begin{cor}
    \label{iotaincl}
    Suppose that $k\ge 1$ and that $p_m\in\{1,\cdots, r-1\}$ for $1\le m\le k $. We have
  $$  W(\iota^-_{p_1}\cdots\iota^-_{p_k}\bos)\hookrightarrow W(\bos)\twoheadrightarrow W(\iota^+_{p_1}\cdots\iota^+_{p_k}\bos).$$
\end{cor}
\subsection{The elements $\bos^\pm$} For $\bos\in \mathbb I_{n}^r$ set $$\bos^\pm=(\iota_{r-1}^\pm\cdots\iota_1^\pm)(\iota_
{r-1}^\pm\cdots\iota_2^\pm)\cdots (\iota^\pm_{r-1}\iota^\pm_{r-2})\iota_{r-1}^\pm\bos.$$
\begin{prop}\label{spm}
 We have $\bos^\pm\in\mathbb I_{n,\pm}^r$.
\end{prop}
\begin{pf}
We show that $\bos^+\in \mathbb I_{n,+}^r$; the proof of $\bos^-\in \mathbb I_{n,-}^r$ is analogous so we omit details.  We proceed by induction on $r$ which holds vacuously true when $r=1$. The inductive hypothesis gives 
$$\bos_1^+ = (\iota_{r-1}^+ \cdots \iota_2^+)\cdots (\iota_{r-2}^+ \iota_{r-1}^+)\iota_{r-1}^+ \bos \ =\  ([i_1,j_1])\vee \bos_2^+,$$
where $\bos_2^+ = ([i_{w(2)},j_{u(2)}], \cdots, [i_{w(r)},j_{u(r)}]) \in \mathbb I_{n,+}^{r-1}$, for some $w,u \in\Sigma_{r-1}$. \\\\ 
If $j_1\geq j_s$, for all $2\leq s\leq r$, it follows that $\bos_1^+\in \mathbb I_{n,+}^r$ and hence $\iota_{r-1}^+\cdots \iota_1^+\bos_1^+ = \bos_1^+$ which completes the proof in this case. Otherwise, since $\bos_2^+\in\mathbb I_{n,+}^{r-1}$ we must have $j_1<j_{u(2)}$ and then
$$\iota_1^+\bos_1^+ = ([\ell_1,j_{u(2)}], [\ell_2,j_1])\vee \bos_1^+(2,r),$$
for a suitable choice of $\ell_1,\ell_2$ such that $\{\ell_1,\ell_2\}=\{i_1, i_{w(2)}\}$. If $j_1\geq j_{u(3)}$ then $\iota_1^+\bos_1 \in \mathbb I_{n,+}^r$ and we are done; otherwise we iterate this process and the proof that $\bos^+\in \mathbb I_{n,+}^r$ is complete. 
\end{pf}

\noindent {\bf Example.} Assume that $n\geq 8$ and let $\bos = ([0,6],[4,8],[2,5])$. Then 
\begin{gather*}
    \iota_2^+ \bos = \bos, \ \ \iota_1^+\iota_2^+ \bos = ([0,8],[4,6],[2,5])  =\iota_2^+\iota_1^+\iota_2^+ \bos = \bos^+ ,\\ 
    \iota_2^-\bos = ([0,6],[4,5],[2,8]), \ \ \iota_1^-\iota_2^-\bos = ([4,5],[0,6],[2,8]) =\iota_2^-\iota_1^-\iota_2^-\bos= \bos^-.
\end{gather*}

\subsection{} The proof of Proposition \ref{mixed} is completed as follows.
It is immediate from Corollary \ref{iotaincl} that we have $$W(\bos^-)\hookrightarrow W(\bos)\twoheadrightarrow W(\bos^+).$$   By Propositions \ref{spm} and \ref{antiweyl} we have an inclusion $V(\bomega_{\bos^-})\hookrightarrow W(\bos^-)$. On the other hand by Proposition \ref{spm} and  Corollary \ref{weylpermute} we have $W(\bos^+)\cong W(\bomega_{\bos^+})$ and hence a map $W(\bos^+)\to V(\bomega_{\bos^+})$ and the proof of Proposition \ref{mixed} is complete.

\begin{rem} In general the socle of a mixed Weyl module is not simple or multiplicity free. The following example is due to Grigorev--Mukhin \cite{GM25} in the case of $\widehat\bu_1$.\\\\
Let $\bos = ([1,2],[2,3],[1,2],[0,1],[1,2],[0,1])$. Then $W(\bos)$ is an example of a mixed Weyl module for  $\widehat\bu_1$ which has a two dimensional family of trivial submodules in its socle. \\\\
In fact, their example can be generalized to arbitrary $n$ as follows. Consider $$\bos=([i,j], [j,n+1+i],[i,j], [-n-1+j,i], [i,j],[-n-1+j, i]).$$
Then $W(\bos)$ contains two copies of the trivial module. The proof follows the same lines as \cite{GM25}; namely one proves that 
\begin{gather*} W(\bos(0,3))\cong \left(V(\bomega_{i,j})\otimes V(\bomega_{\bos(1,3)})\right)\oplus V(\bomega_{i,j})\cong V(\bomega_{\bos(0,3)})\oplus V(\bomega_{i,j}),\\
W(\bos(3,6))\cong V(\bomega_{\bos(3,6)})\oplus V(\bomega_{-n-1+j,i}).
\end{gather*}
Noting that $$V(\bomega_{-n-1+j,i})^* \cong V(\bomega_{i,j}) \ \ {\rm and} \ \ V(\bomega_{\bos(3,6)})^*\cong V(\bomega_{\bos(0,3)}),$$ we see that the trivial occurs twice in $W(\bos)\cong W(\bos(0,3))\otimes W(\bos(3,6))$.
\end{rem}

\subsection{Extensions}\label{extensions}
As another application of our results we discuss extensions in the category $\mathscr F_n$. Given a closed subset $\mathbb J\subset\mathbb I_n^r$, let $\mathscr F_n(\mathbb J)$ be the full subcategory of $\mathscr F_n$ consisting of objects whose Jordan--H\"older constituents are isomorphic to $V(\bomega_{\bos})$ with $\bomega_\bos\in\{\bomega_{\bos'}:\bos'\in\mathbb J\}$. Clearly $\mathscr F_n$ is an abelian category.
\begin{prop}
   Suppose that $\mathbb J_1$ and $\mathbb J_2$ are closed subsets of $\mathbb I_{n}^{r_1}$ and $ \mathbb I_{n}^{r_2}$, respectively, and $$\{\bomega_\bos:\bos\in\mathbb J_1\}\cap\{\bomega_\bos:\bos\in\mathbb J_2\}=\emptyset.$$ Suppose that $V_p$ is an object of $\mathscr F_n(\mathbb J_p)$ for $p=1,2$. Then
   $$\Ext^s_{\widehat\bu_n}(V_1, V_2)=0,\ \  s=0,1.$$ 
\end{prop}
\begin{pf} It is clear from Theorem \ref{weylmaps} that $\Hom_{\widehat\bu_n}(V_1,V_2)=0$ and hence the proposition holds for $s=0$.
\\\\
For $s=1$ we must  prove that any short exact sequence $$0\to V_2\to V\to V_1\to 0$$ of modules in $\mathscr F_n$ splits.
    Assume first that $V_p= V(\bomega_{\bos_p})$ are objects of $\mathscr F_n(\mathbb J_p)$ for $p=1,2$. If $V$ is a quotient  of $W(\bomega_\bos)$ then since $V(\bomega_{\bos_1})$ is a further quotient of $V$, we must have $\bomega_\bos=\bomega_{\bos_1}$. Since $V(\bomega_{\bos_2})$ is a submodule of $V$ we have $\bomega_{\bos_2}\in\wt_\ell^+ W(\bomega_{\bos_1})$. It follows from Theorem \ref{weylmaps} that $\bomega_{\bos_2}=\bomega_{\bos'}$ for some $\bos'\in\bar\bos_1[n]$. Hence we get 
    $$\{\bomega_\bos:\bos\in\mathbb J_2\}\ni\bomega_{\bos_2}=\bomega_{\bos'}\in\{\bomega_{\bos''}:\bos''\in\bar\bos_1[n]\}\subset\{\bomega_\bos:\bos\in\mathbb J_1\},$$ where the last inclusion follows since $\mathbb J_1$ is closed. But this contradicts the hypothesis of the proposition.
    \\\\
    In particular, this proves that  $\bomega_{\bos_2}\notin\bomega_{\bos_1}(\cal Q^+_n)^{-1}$. Working with duals, a similar argument proves that $\bomega_{\bos_1}\notin\bomega_{\bos_2}(\cal Q_n^+)^{-1}$. Hence it follows that the pre-image of $V(\bomega_{\bos_1})$ cannot intersect $V(\bomega_{\bos_2})$ and the proof is complete in this case.
    \\\\
    Suppose now that we know the result when $V_1$ has length $\ell-1$. To prove it when it has length  $\ell$, let $V(\bomega_{\bos_1})$ be a quotient of $V_1$ and  consider the short exact sequence $$0\to V_1'\to V_1\to V(\bomega_{\bos_1})\to 0.$$ Since $\Hom_{\widehat\bu_n}(V_1, V(\bomega_{\bos_2}))=0$ we have a left exact sequence
    $$0\to \Ext^1(V(\bomega_{\bos_1}), V(\bomega_{\bos_2}))\to \Ext^1(V_1, V(\bomega_{\bos_2}))\to \Ext^1(V_1', V(\bomega_{\bos_2})).$$
    Since the terms on the ends are zero by the inductive hypothesis, it follows that the middle term is zero as needed. A similar argument applies if we replace $V(\bomega_{\bos_2})$ by $V_2$ and the proposition follows.
\end{pf}

\providecommand{\bysame}{\leavevmode\hbox to3em{\hrulefill}\thinspace}
\providecommand{\MR}{\relax\ifhmode\unskip\space\fi MR }
\providecommand{\MRhref}[2]{%
  \href{http://www.ams.org/mathscinet-getitem?mr=#1}{#2}
}
\providecommand{\href}[2]{#2}


\begin{thebibliography}{10}

\bibitem{BC25}
Matheus Brito and Vyjayanthi Chari, \emph{Alternating snake modules and a determinantal formula}, arXiv:2412.03750.

\bibitem{C01}
Vyjayanthi Chari, \emph{On the fermionic formula and the {K}irillov-{R}eshetikhin conjecture}, Internat. Math. Res. Notices (2001), no.~12, 629--654. \MR{1836791 (2002i:17019)}

\bibitem{Ch01}
\bysame, \emph{Braid group actions and tensor products}, Int. Math. Res. Not. (2002), no.~7, 357--382. \MR{1883181}

\bibitem{CL06}
Vyjayanthi Chari and Sergei Loktev, \emph{Weyl, {D}emazure and fusion modules for the current algebra of {$\mathfrak{sl}_{r+1}$}}, Adv. Math. \textbf{207} (2006), no.~2, 928--960. \MR{2271991 (2008a:17029)}

\bibitem{CP91}
Vyjayanthi Chari and Andrew Pressley, \emph{Quantum affine algebras}, Comm. Math. Phys. \textbf{142} (1991), no.~2, 261--283. \MR{1137064}

\bibitem{CP95}
\bysame, \emph{Quantum affine algebras and their representations}, Representations of groups ({B}anff, {AB}, 1994), CMS Conf. Proc., vol.~16, Amer. Math. Soc., Providence, RI, 1995, pp.~59--78. \MR{1357195 (96j:17009)}

\bibitem{CP01}
\bysame, \emph{Weyl modules for classical and quantum affine algebras}, Represent. Theory \textbf{5} (2001), 191--223 (electronic). \MR{1850556 (2002g:17027)}

\bibitem{FoL07}
Ghislain Fourier and Peter Littelmann, \emph{Weyl modules, {D}emazure modules, {KR}-modules, crystals, fusion products and limit constructions}, Adv. Math. \textbf{211} (2007), no.~2, 566--593. \MR{2323538 (2008k:17005)}

\bibitem{FM01}
Edward Frenkel and Evgeny Mukhin, \emph{Combinatorics of q-characters of finite-dimensional representations of quantum affine algebras}, Communications in Mathematical Physics \textbf{216} (2001), no.~1, 23--57.

\bibitem{FR99}
Edward Frenkel and Nicolai Reshetikhin, \emph{The {$q$}-characters of representations of quantum affine algebras and deformations of {$\mathcal{W}$}-algebras}, Recent developments in quantum affine algebras and related topics ({R}aleigh, {NC}, 1998), Contemp. Math., vol. 248, Amer. Math. Soc., Providence, RI, 1999, pp.~163--205. \MR{1745260}

\bibitem{GM25}
Andrei Grigorev and Evgeny Mukhin, \emph{On representations of quantum affine ${\mathfrak{sl}}_2$}, in preparation.

\bibitem{Gur21}
Maxim Gurevich, \emph{On the hecke-algebraic approach for general linear groups over a $p$-adic field}, Interactions of Quantum Affine Algebras with Cluster Algebras, Current Algebras and Categorification, in honor of Vyjayanthi Chari, Progress in Mathematics 337, 2021.

\bibitem{HL10}
David Hernandez and Bernard Leclerc, \emph{Cluster algebras and quantum affine algebras}, Duke Math. J. \textbf{154} (2010), no.~2, 265--341. \MR{2682185}

\bibitem{HL13a}
\bysame, \emph{A cluster algebra approach to {$q$}-characters of {K}irillov-{R}eshetikhin modules}, Journal of the European Mathematical Society \textbf{18} (2013).

\bibitem{KKOP}
Masaki Kashiwara, Myungho Kim, Se-jin Oh, and Euiyong Park, \emph{Monoidal categorification and quantum affine algebras}, Compositio Mathematica \textbf{156} (2020), no.~5, 1039–1077.

\bibitem{KKOP22}
Masaki Kashiwara, Myungho Kim, {Se jin} Oh, and Euiyong Park, \emph{Cluster algebra structures on module categories over quantum affine algebras}, Proceedings of the London Mathematical Society \textbf{124} (2022), no.~3, 301--372 (English).

\bibitem{LM14}
Erez Lapid and Alberto Mínguez, \emph{On a determinantal formula of {T}adi\'c}, American Journal of Mathematics \textbf{136} (2014), 111--142.

\bibitem{LM18}
\bysame, \emph{Geometric conditions for $\square$-irreducibility of certain representations of the general linear group over a non-archimedean local field}, Advances in Mathematics \textbf{339} (2018), 113--190.

\bibitem{MY12}
Evgeny Mukhin and Charles A.~S. Young, \emph{Path description of type ${B}$ $q$-characters}, Advances in Mathematics \textbf{231} (2012), 1119--1150.

\bibitem{VV02}
M.~Varagnolo and E.~Vasserot, \emph{Standard modules of quantum affine algebras}, Duke Math. J. \textbf{111} (2002), no.~3, 509--533. \MR{1885830}

\end{thebibliography}
\end{document}